\documentclass[MSNbibl,number,citesort,seceqn,dvips]{arxbj}
\usepackage{upgreek}

\aid{0}
\volume{20}
\issue{3}
\pubyear{2014}
\firstpage{1234}
\lastpage{1259}
\doi{10.3150/13-BEJ520} %kopijuoti is PTS

\makeatletter
\newcommand{\RMo}{\mathrm{o}}
\newcommand{\RMO}{\mathrm{O}}

\newcommand{\RMi}{\mathrm{i}}

\newcommand{\mrmd}{\,\mathrm{d}}

\newproclaim{assumption}{Assumption}

\newtheorem{claim}{Claim}

\newcommand{\lldots}{\ldots} %%%% neaisku
\newcommand{\rrvert}{\vert}
\newcommand{\llvert}{\vert}

\newcommand{\cal}{\mathcal}

\def\E{\mathbb{E}}
\newcommand{\Tr}{\operatorname{Tr}}

\newtheorem{lem}{Lemma}
\newtheorem{theorem}{Theorem}[section]

\newremark{rem}{Remark}[section]

\newcommand{\I}{\mathbb{I}}

\makeatother

\begin{document}
\begin{frontmatter}

\title{Limiting spectral distribution of sample autocovariance matrices}
\runtitle{Autocovariance matrix}

\begin{aug}
%%%% inicialai - be tarpu
\author[1]{\inits{A.}\fnms{Anirban} \snm{Basak}\thanksref{1}\ead[label=e1]{anirbanb@stanford.edu}}, % \and
\author[2]{\inits{A.}\fnms{Arup} \snm{Bose}\corref{}\thanksref{2}\ead[label=e2]{bosearu@gmail.com}} \and
\author[3]{\inits{S.}\fnms{Sanchayan} \snm{Sen}\thanksref{3}\ead[label=e3]{sen@cims.nyu.edu}}
\runauthor{A. Basak, A. Bose and S. Sen} %% auto
\address[1]{Department of Statistics, Stanford University, 390 Serra Mall, Stanford, CA 94305-4065, USA.\\
\printead{e1}}
\address[2]{Statistics and Mathematics Unit, Indian Statistical
Institute, 203 B. T. Road, Kolkata 700108, India. \printead{e2}}
\address[3]{Courant Institute of Mathematical Sciences, New York University, 251 Mercer Street, New
York, NY 10012, USA. \printead{e3}}
\end{aug}

% HISTORY:
\received{\smonth{9} \syear{2011}}
\revised{\smonth{12} \syear{2012}}

% ABSTRACT
%
\begin{abstract}
We show that the empirical spectral distribution (ESD) of the sample
autocovariance matrix (ACVM) converges as the dimension increases, when
the time series is a linear process with reasonable restriction on the
coefficients. %This
The limit does not depend on the distribution of the underlying driving
i.i.d. sequence and its support is unbounded. This limit does not
coincide with the spectral distribution of the theoretical ACVM.
However, it does so if we consider a suitably tapered version of the
sample ACVM. For banded sample ACVM the limit has unbounded support as
long as the number of non-zero diagonals in proportion to the dimension
of the matrix is bounded away from zero. If this ratio tends to zero,
then the limit exists and again coincides with the spectral
distribution of the theoretical ACVM. Finally, we also study the LSD of
a naturally modified version of the ACVM which is not non-negative
definite.
\end{abstract}

% KEYWORDS
% visi is mazosios raides ir pagal abecele
%
\begin{keyword}
\kwd{autocovariance function}
\kwd{autocovariance matrix}
\kwd{banded and tapered autocovariance matrix}
\kwd{linear process}
\kwd{spectral distribution}
\kwd{stationary process}
\kwd{Toeplitz matrix}
\end{keyword}

\end{frontmatter}

%s1 #&#
\section{Introduction}\label{sectionIntroduction}
Let $X=\{X_t\}$ be a \textit{stationary} process with $\E(X_t) = 0$
and $\E(X_t^2) < \infty$. The \textit{autocovariance function}
(ACVF) $\gamma_X(\cdot)$ and the \textit{autocovariance matrix}
(ACVM) $\Sigma_n(X)$ of order $n$ are defined as:
\[
\gamma_X(k) = \operatorname{cov}(X_0, X_k),\qquad k =0, 1,
\ldots
\]
and
\[
\Sigma_n (X) = \bigl(\bigl(
\gamma_X(i-j)\bigr)\bigr)_{1 \leq i,j \leq n}.
\]
To every ACVF,
there corresponds a unique distribution, called the \textit{spectral
distribution}, $F_X(\cdot)$ which satisfies
%
%e1.1 #&#
\begin{equation}
\label{eqspectraldistribution}\gamma_X(h)=\int_{(0,1]}
\exp(2\uppi  \RMi  hx)\mrmd F_X(x) \qquad\mbox{for all } h.
\end{equation}
We shall assume that
%
%e1.2 #&#
\begin{equation}
\label{eqsummability}\sum_{k=1}^\infty\bigl|
\gamma_X(k)\bigr| < \infty.
\end{equation}
Then $F_X(\cdot)$ has a density,
known as
the \textit{spectral density of} $X$ or of $\gamma_X(\cdot)$, which equals
%
%e1.3 #&#
\begin{equation}
\label{eqspectraldensity}f_X(t)=\sum_{k=-\infty
}^\infty
\exp{(-2\uppi \RMi tk)}\gamma_X(k),\qquad t\in(0, 1].
\end{equation}
The
\textit{non-negative definite} estimate of $\Sigma_n(X)$ is the
\textit{sample ACVM}
%
%e1.4 #&#
\begin{equation}
{\Gamma}_n (X)= \bigl(\bigl(\hat{\gamma}_X(i-j)\bigr)
\bigr)_{1 \leq i,j \leq n} \qquad\mbox{where } \hat{\gamma}_X(k)=
n^{-1} \sum_{i=1}^{n-|k|}
X_{i}X_{i+|k|}.
\end{equation}
%
%When $\{\varepsilon_t\}$ is i.i.d.,

The matrix $\Gamma_n(X)$ is a random matrix. Study of the behavior of
random matrices, when the dimension goes to $\infty$,
have been inspired by both theory and applications. This is done by
studying the behavior of its eigenvalues. For instance a host of
results are known for the related sample covariance matrix, in the
i.i.d. set-up and its variations; results on its spectral distribution,
spacings of the eigenvalues, spectral statistics etc. encompasses a
rich theory and a variety of applications.

The autocovariances are of course crucial objects in time series analysis.
They are used in estimation, prediction, model fitting and white noise tests.
Under suitable assumptions on $\{X_t\}$, for every fixed $k$,
$\hat{\gamma}_X(k) \to\gamma_X(k)$ almost surely (a.s.). There
are also results on the asymptotic distribution of specific functionals
of the autocovariances. Recently, there has been growing interest in
the matrix $\Gamma_n(X)$ itself.
For instance,
the largest
eigenvalue
of $\Sigma_n(X)-\Gamma_n(X)$
does not converge to
zero, even under reasonable assumptions (see Wu and Pourahmadi \cite
{WuPourahmadi}, Arcones \cite{mcmurraypolitis2010} and Xiao and Wu
\cite{xiaowu}).

In this article we study the behavior of $\Gamma_n(X)$, and a few
other natural estimators of $\Sigma_n(X)$, as $n \to
\infty$, through the behavior of its spectral
distribution. We investigate the \emph{consistency} (in an appropriate
sense) of these estimators.
%for the \emph{limiting spectral distribution}
%for $\Sigma_n(x)$.

For a real symmetric matrix $A_{n \times n}$ with eigenvalues $\lambda
_1$, $\lambda_2,\ldots,\lambda_n$ the \textit{Empirical Spectral
Distribution} (\textit{ESD}) of $A_n$ is
defined as,
%
%e1.5 #&#
\begin{equation}
F^{A_{n}}(x)= n^{-1}\sum_{i=1}^{n}
\I(\lambda_i\le x).
\end{equation}
If $\{F^{A_n}\}$ converges weakly to $F$, we
write $F^{A_n} \stackrel{w}{\to} F$. For $X$ any random variable with
distribution $F$, $X$ or $F$ will be called the \textit{Limiting
Spectral Distribution} (\textit{or measure}) (\textit{LSD}) of $F^{A_n}$. The entries of
$A_n$ are allowed to be random. In that case, the limit is taken to be either
in probability or (as in this paper) in a.s. sense.

Any matrix $T_n$ of the form
$((t_{i-j}))_{1\leq i, j \leq n}$ is a \textit{Toeplitz} matrix and
hence $\Sigma_n(X)$
and $\Gamma_n(X)$ (with a triangular sequence of entries) are
Toeplitz matrices.
For $T_n$ symmetric, from
Szeg{\"o}'s theory of Toeplitz operators
(see B{\"o}ttcher and Silbermann \cite{bottchersilberman}), we note that
if $\sum|t_k | < \infty$, then the LSD
of $T_n$
equals
$f(U)$ where $U$ is uniformly distributed on $(0, 1]$ and
$f(x)= \sum_{k=-\infty}^{\infty} t_k \exp{(-2\uppi  \RMi  x k)}$, $ x\in
(0, 1]$.
In particular if (\ref{eqsummability}) holds, then the LSD of
$\Sigma_n(X)$ equals
$f_X(U)$ where $f_X(\cdot)$ is as defined in (\ref{eqspectraldensity}).

We call a sequence of estimators $\{E_n\}$ of $\Sigma_n(X)$
\textit{consistent} if its LSD is $f_X(U)$ where $U$ is uniformly
distributed on $[0, 1]$. We show that $\{\Gamma_n(X)\}$ is inconsistent
(see Theorem \ref{thmautocov}(c)). We also show that if $\Gamma_n(X)$
is modified by suitable tapering or banding then the modified
estimators are indeed consistent (see Theorem \ref {thmautocovband}(b)
and (c)). This phenomenon is mainly due to the estimation of a large
number of autocovariances by $\Gamma_n(X)$. Such
inconsistency of sample covariance matrices has also been observed in
the context of high-dimensional multivariate analysis, and is now well
understood, with the help the results from Random Matrix Theory.

To obtain the convergence of ESD of such estimators, we impose a
reasonable condition on the stationary process $\{X_t\}$; we assume it
to be a linear process, that is,
%
%e1.6 #&#
\begin{equation}
\label{eqsetup}X_t=\sum_{k=0}^{\infty}
\theta_k\varepsilon_{t-k},
\end{equation}
where $\{\theta_k\}$ satisfies a weak condition and $\{\varepsilon_t,
t\in\mathbb{Z}\}$ is a sequence of
independent random variables with appropriate conditions.
The simulations of Sen \cite{SenA} suggested that the LSD of $\Gamma_n(X)$
exists and is independent of
the distribution of $\{\varepsilon_t\}$ as long as they are
i.i.d. with mean zero and variance one.
Basak \cite{Basak} and Sen \cite{SenS} initially studied,
respectively, the special cases where $X$ is an i.i.d. process or is
an MA(1) process.

In Theorem \ref{thmautocov}, we prove that, if
$\{X_t\}$ satisfies (\ref{eqsetup}) and $\sum_{k=0}^\infty|\theta
_k| < \infty$ then the LSD of $\Gamma_n(X)$ exists, and it is \textit
{universal} when
$\{\varepsilon_t\}$ are independent with mean zero and variance 1
and are either uniformly bounded or identically distributed. We further
show that LSD is unbounded when $\theta_i \geq0$ for all $i$, and
thus $\{\Gamma_n(X)\}$ is inconsistent, since $f_X(U)$ is of bounded support.

When $\{X_t\}$ is a finite order process, the limit moments
%in our case
can be
written as multinomial type sums of the
autocovariances (see (\ref{eqnautocov-limit-moment})).
When $X$ is of infinite order, the limit moments are the limits
%values
of
these sums
%multinomial expressions
as the order tends to
infinity. Additional properties of the limit moments are
available in the companion report Basak, Bose and Sen \cite
{Basakbosesen}.

Incidentally, $\Gamma_n(X)$ reminds us of the
sample
covariance matrix, $S$, for the i.i.d. set-up, whose spectral properties
are well known. See
Bai \cite{Bai99} for
%some of
the basic
references on $S$. In particular, the LSD of $S$ (with i.i.d. entries) under
suitable conditions
is the Mar\v{c}enko--Pastur law and
%given by Mar{\v{c}}enko and Pastur \cite{marc}, which
is supported on the interval $[0, 4]$. Thus,
%our result on
the LSD of $\Gamma_n(X)$ is in sharp contrast.

The proof of Theorem \ref{thmautocov} is challenging, mainly because
of the non-linear dependence, and the Teoplitz structure of $\Gamma
_n(X)$. Bai and Zhou \cite{Baizhou2008} and Yao \cite{Yao2012} study
the LSD of the sample covariance matrix of $\mathbf{X}_1,\ldots,
\mathbf{X}_n$ where $\mathbf{X}_k$ are i.i.d. $p$-dimensional vectors with
some dependence structure. They establish the existence of the LSD by
using Stieltjes transform method. Here this approach fails completely
due to the strong row column dependence. In fact no Stieltjes transform
proof for even the Toeplitz matrix with i.i.d. input is known.
Moreover one added advantage in both the above articles is the
existence of $n$ independent columns, which we lack here, because we
have only one sample from the linear process $\{X_t\}$. The methods of
Xiao and Wu \cite{xiaowu} is also not applicable in our set-up because
they deal with only the maximum eigenvalue of the difference of $\Sigma
_n(X)$, and $\Gamma_n(X)$, not the ESD of $\Gamma_n(X)$.

%The
%Note that the
%LSD of
%$\Sigma_n(X)$ depends on the parameters $\{\theta_k\}$ but
%there is no one to one correspondence between them.
%$\{\theta_k\}$ and
%the LSD--
%For instance,
%the LSD is same when $X$ is AR(1) with
%parameter $\theta$ or $-\theta$. The same situation persists for the
%LSD of $\Gamma_n(X)$ (see Remark \ref{remkernel}(iii)).

%We call a sequence of estimators $\{E_n\}$ of $\Sigma_n(X)$ to be
%distributed on $[0, 1]$. Thus $\{\Gamma_n(X)\}$ is inconsistent.
Now consider a sequence of integers $m:=m_n\to\infty$, and a kernel
function $K(\cdot)$. Define
%
%e1.7 #&#
\begin{equation}
\label{eqkde}\hat f_X(t)=\sum_{k=-m}^m
K(k/m)\exp{(-2\uppi \RMi tk)}\hat\gamma_X(k),\qquad t\in(0, 1]
\end{equation}
as the kernel density estimate of $f_X(\cdot)$.
Considering this as a spectral
density, the corresponding ACVF
is given by (for $-m \leq h \leq m$):
\begin{eqnarray*}
\gamma_K(h)&=&\int_{(0, 1]} \exp(2\uppi  \RMi  hx) \hat
{f}_X(x)\mrmd x
\\
&=& \sum_{k=-m}^m K(k/m)\int
_{(0, 1]} \exp\{2\uppi  \RMi  hx-2\uppi \RMi xk\}\hat\gamma_X(k)\mrmd x\\
&=&
K(h/m) \hat\gamma_X(h) % \mbox{for all} -m \leq j \leq m.
\end{eqnarray*}
and is $0$ otherwise. This motivates the consideration of the
\textit{tapered sample ACVM}
%
%e1.8 #&#
\begin{equation}
\label{eqtaperedauto}\Gamma_{n,K}(X)=\bigl(\bigl(K\bigl((i-j)/m\bigr)
\hat\gamma_X(i-j)\bigr)\bigr)_{1\leq i, j\leq n}.
\end{equation}
%
%It may be noted that
If $K$ is a non-negative definite function
then $\Gamma_{n,K}(X)$ is also non-negative definite.
Among other results, Xiao and Wu \cite{xiaowu} also showed that under the
growth condition $m_n=\RMo(n^\gamma)$ for a suitable $\gamma$
and suitable conditions on $K$, the largest eigenvalue of
$\Gamma_{n,K}(X)-\Sigma_n(X)$ tends to zero a.s. Theorem \ref
{thmautocovband}(c) states that under the minimal condition $m_n/n
\to0$, if $K$ is bounded, symmetric and continuous at 0 and
$K(0)=1$, then $\Gamma_{n, K}(X)$ is consistent.
This is a reflection of the fact that the consistency notion of Xiao
and Wu \cite{xiaowu} in terms of the maximum eigenvalue is stronger
than our notion and hence our consistency holds under weaker growth
condition on $m_n$.

The second approach is to use banding as in McMurry and Politis \cite
{mcmurraypolitis2010} %use such banded matrices
who used it to
develop their bootstrap
procedures. We study two such banded matrices. Let $\{m_n\}_{n \in
\mathbb{N}}\to\infty$ be such that $\alpha_n:= m_n/ n \to\alpha
\in[0, 1]$. Then the \textit{type I banded sample autocovariance
matrix} $\Gamma_n^{\alpha, I}(X)$ is same as $\Gamma_n(X)$ except
that we substitute $0$ for $\hat{\gamma}_X(k)$ whenever $|k| \ge
m_n$. This is the same as
$\Gamma_{n,K}$ with
$K(x)=I_{\{|x|\leq1\}}$. The \textit{type II banded
ACVM}
$\Gamma_n^{\alpha, \mathit{II}}(X)$ is the $m_n \times
m_n$ principal sub matrix of $\Gamma_n(X)$. Theorem
\ref{thmautocovband}(a) and (b) states our results on these banded
ACVMs.
In particular, the LSD exists for all
$\alpha$ and is unbounded when $\alpha\ne0$.
When $\alpha=0$, the LSD is $f_X(U)$ and thus those estimate matrices
are consistent.

A related matrix, which may be of
%some
interest, especially to probabilists,
is,
%
%e1.9 #&#
\begin{equation}
\Gamma_{n}^{*}(X) = \bigl(\bigl( \gamma^{*}_X\bigl(|i-j|\bigr)
\bigr)\bigr)_{1 \leq i,j \leq n } \qquad\mbox{where } \gamma_X^{*}(k)=
n^{-1} \sum_{i=1}^{n}
X_{i} X_{i+k}, k =0, 1, \ldots.
\end{equation}
$\Gamma_n^*(X)$ does not have a ``data'' interpretation unless one
assumes we have $2n-1$ observations $X_1,\ldots, X_{2n-1}$. It is not
non-negative definite and hence
many of the techniques applied to $\Gamma_n(X)$ are not available for
it.
Theorem \ref{thmautocov2} states that its LSD also exists but under
stricter conditions on $\{X_t\}$. Its moments dominate those of the LSD
of $\Gamma_n(X)$ when $\theta_i \geq0$ for all $i$ (see Theorem~\ref
{thmautocov2}(c))
even though simulations show that the LSD of $\Gamma_n^*(X)$ has
significant positive mass on the negative axis. %
%To illustrate our results,
%%in Section \ref{subsectionsimulations}
%we provide a few simulation results for different choices of
%$\{\theta_k\}$. It would be nice to obtain additional theoretical
%properties of the ESD and the LSD of these matrices. For instance,
%the distribution of maximum eigenvalue of $S$ has been studied in
%the literature. However, it does not seem to be at all easy to
%obtain similar results for $\Gamma_n(X)$.

%developed in Section
%Apart from providing more information on the nature
%of the limit, some of these results are used crucially in the proof
%of Theorem \ref{thmlsdnonunique}.

%s2 #&#
\section{Main results}\label{sectionMainresults}
We shall assume that $X=\{X_t\}_{t\in\mathbb{Z}}$ is a linear
(MA($\infty$)) process
%
%e2.1 #&#
\begin{equation}
\label{deflinear} X_t=\sum_{k=0}^{\infty}
\theta_k\varepsilon_{t-k},
\end{equation}
where $\{\varepsilon_t, t\in\mathbb{Z}\}$ is a sequence of
independent random variables.
A special case of this process is the so called
MA($d$) where $\theta_k=0$ for all $k >d$. We denote this process by
\[
X^{(d)}=\{X_{t,d}\equiv\theta_{0}
\varepsilon_{t}+\theta_{1}\varepsilon_{t-1} +
\cdots+ \theta_{d}\varepsilon_{t-d}, t \in\mathbb{Z}\}\qquad
(\theta_0\neq0).
\]
%
%where without loss we assume that $$.
Note that
%It may also be mentioned that
working with two sided moving average entails no difference. The
conditions on $\{\varepsilon_t\}$ and on $\{\theta_k\}$ that will be
used are:

\renewcommand{\theassumption}{\Alph{assumption}}
\begin{assumption}\label{assumptionA}
\textup{(a)} $\{\varepsilon_t\}$ are i.i.d. with $\E[\varepsilon
_{t}]=0$ and $\E[\varepsilon_{t}^2]=1$.

\textup{(b)} $\{\varepsilon_t\}$ are independent, uniformly
bounded with $\E[\varepsilon_{t}]=0$ and
$\E[\varepsilon_{t}^2]=1$.
\end{assumption}

\begin{assumption}\label{assumptionB}
\textup{(a)} $\theta_j\geq0$ for all $j$.
%(b) $\theta_j\geq0$ for all $0 \leq j < \infty$.

\textup{(b)} $\sum_{j=0}^\infty|\theta_j|<\infty$.
\end{assumption}

%It may be noted that
The series in (\ref{deflinear}) converges
a.s. under Assumptions \ref{assumptionA}(a) (or (b)) and
%Assumption
\ref{assumptionB}(b).
Further, $X$ and $X^{(d)}$ are strongly stationary and ergodic
under Assumption \ref{assumptionA}(a) and weakly (second order) stationary under
Assumptions \ref{assumptionA}(b) and
%Assumption
\ref{assumptionB}(b).

The ACVF of
$X^{(d)}$ and $X$ are given by
%
%e2.2 #&#
\begin{equation}
\label{eqndef-gammaj}\gamma_{X^{(d)}}(j)=\sum
_{k=0}^{d-j}\theta_k
\theta_{j+k} \quad\mbox{and}\quad \gamma_{X}(j)=\sum
_{k=0}^{\infty}\theta_k
\theta_{j+k}.
\end{equation}
Let $\{k_i\}$ stand for suitable integers and let
%
%e2.3 #&#
\begin{equation}
\label{eqkShd}\mathbf{k}=(k_0,\lldots, k_d),\qquad
S_{h,
d}=\{\mathbf{k}\dvt  k_0,\ldots, k_d\geq0,
k_0+ \cdots+ k_d=h\}.
\end{equation}
%
%Suppose Assumption A(a) or A(b) holds and $ 0 \leq\alpha\leq1$.
%(a) Then almost surely, $F^{\Gamma_{n}^\alpha(X^{(d)})} \stackrel{w}{
%depend on the distribution of $\{\varepsilon_{t}\}$.
%Further, for some sequence of constants $\{p^{\alpha, (d)}_{
%(d)}_{\mathbf{k}}\prod_{i=0}^d[\gamma_{X^{(d)}}(i)]^{k_i}.
%(b) Under Assumption B(c), almost surely, $F^{\Gamma{^\alpha_{n}}(X)}
%the distribution of $\{\varepsilon_{t}\}$. Further as $d\to\infty$,
%$$F^\alpha_d\stackrel{w}{\to}F^\alpha\mbox{and} \beta^
%(c) Under Assumption B(a), if $ \alpha> 0$, then $F^\alpha_d$ has
%unbounded support and $\beta^\alpha_{h,d-1}\leq\beta^\alpha_{h,d}$
% if $d\geq1$. Hence if
%Assumption B(b) and B(c) hold, then $F^\alpha$ has unbounded
%support for $\alpha> 0$.
%If $\alpha=0$, then all the limits above have bounded support.

%th2.1 #&#
\begin{theorem}[(Sample ACVM)]\label{thmautocov}
Suppose Assumption \textup{\ref{assumptionA}(a)} or \textup{(b)} holds. %$ 0 \leq\alpha\leq1$.

\textup{(a)} Then a.s., $F^{\Gamma_{n}(X^{(d)})} \stackrel{w}{\to}F_d$ which
is non-random and does not depend on the distribution of $\{\varepsilon
_{t}\}$.
Further,
%
%e2.4 #&#
\begin{equation}
\label{eqnautocov-limit-moment} \beta_{h,d}=\int x^h \mrmd
F_d(x)=\sum_{S_{h,d}}p^{(d)}_{\mathbf
{k}}
\prod_{i=0}^d\bigl[\gamma_{X^{(d)}}(i)
\bigr]^{k_i},
\end{equation}
where $\{p^{(d)}_{\mathbf{k}}\}$ are universal constants independent
of the $\theta_i$ and the $\{\epsilon_i\}$. They are defined by a
limiting process given in (\ref{eqdefinepw}) and (\ref{eqlimitpw}).

\textup{(b)} Under Assumption \textup{\ref{assumptionB}(b)}, a.s., $F^{\Gamma_{n}(X)}\stackrel{w}{\to}
F$ which is non-random and independent of the distribution of $\{
\varepsilon_{t}\}$. Further for every fixed $h$, as $d\to\infty$,
\[
F_d\stackrel{w} {\to}F \quad\mbox{and}\quad \beta_{h,d}
\rightarrow\beta_h=\int{x^h \mrmd F(x)}.
\]

\textup{(c)} Under Assumption \textup{\ref{assumptionB}(a)}, $F_d$ has unbounded support and $\beta
_{h,d-1}\leq\beta_{h,d}$ if $d\geq1$.
%As a
Consequently,
%consequence,
if Assumption \textup{\ref{assumptionB}(a)} and \textup{(b)} holds,
then $F$ has unbounded support. Therefore $\{\Gamma_n(X)\}$ is
inconsistent.
\end{theorem}
%
%We now state an LSD result for $\Gamma_n^*(X)$.
%As mentioned before, $\Gamma_n^*(X)$
%is not non-negative definite and this creates technical difficulties.
%We have the following restrictive result.
%deal with only the case when Assumption A(b) holds
%the $\{\varepsilon_i\}$ are uniformly bounded. Also,
%and
%For simplicity we
%further
%assume that $\alpha=1$.
%
%th2.2 #&#
\begin{theorem}\label{thmautocov2}
Suppose Assumption \textup{\ref{assumptionA}(b)} holds. Then conclusions of Theorem \ref
{thmautocov} continue to hold for $\Gamma_n^*(X)$, $d \leq\infty$,
and (\ref{eqnautocov-limit-moment}) holds with
modified universal constants $\{p_{\mathbf{k}}^{*(d)} \}$.
%
%(a) Then almost surely, $F^{\Gamma_{n}^*(X^{(d)})}
%depend on the distribution of $\{\varepsilon_{t}\}$. For some
%constants $\{p_{\mathbf{k}}^{\ast(d)}\}$,
%(b) Under Assumption B(b), almost surely
%$F^{\Gamma_{n}^*(X)}\stackrel{w}{\to} F^{\ast}$ which is also
%nonrandom and does not depend on the distribution of
%$\{\varepsilon_{t}\}$. Further as $d\to\infty$,
%$$F_{d}^{\ast}\stackrel{w}{\to}F^{\ast} \mbox{and}
%dF^{\ast}(x).$$ (c) Under Assumption B(a),
%$F_{d}^{\ast}$ has unbounded support and $\beta_{h,d-1}^\ast\leq
%Assumption B(a) and B(b), $F^\ast$ has unbounded support. Moreover
%$\beta_h\leq\beta_h^*$ for all $h$.
\end{theorem}

%re2.1 #&#
\begin{rem} \label{remthintails}
\textup{(i)}
From the proofs,
%of the above theorems,
it will follow that the limit moments $\{\beta_{h,d}\}$ and $\{\beta
_{h}\}$
of the above LSDs are dominated by
%$\{p_{\mathbf{k}}^{(d)}\}$ of Theorem \ref{thmautocov} satisfies
%$$p^{(d)}_{\mathbf{k}}\leq
%As a consequence,
$
%|\beta_{h,d}| \leq
\frac{4^h(2h)!}{h!}(\sum_{k=0}^\infty|\theta_k|)^{2h}$
% \mbox{and} |\beta_{h}| \leq\frac{4^h(2h)!}{h!}(\sum_{k=0}^
which are the $(2h)$th moment of a Gaussian variable with mean zero and variance
$4((\sum_{k=0}^\infty|\theta_k|)^2)$.
% of a Gaussian random variable.
%Hence the limits have subexponential tails.
%The same is true for the LSD of $\Gamma_n^*(X)$.
%Moreover
Hence the limit moments uniquely identify the LSDs.
% obtained via the moment convergence are uniquely identified.

\mbox{}\hphantom{i}\textup{(ii)}
%Observe that the expressions for moments in
All the above LSDs have
%The distributions given by the moments as in (
%(\ref{eqnautocov2-limit-moment})
unbounded support while
%the random variable
%are similar to the moments of
%$f_X(U)$ given in (\ref{eqfxumoments}). However, while the former
%two variables have unbounded support
%The variable $f_X(U)$
%latter
$f_X(U)$
has support
contained in $[-\sum_{-\infty}^\infty|\gamma_X(k)|,
\sum_{-\infty}^\infty|\gamma_X(k)|]$.
Simulations show that the LSD of $\Gamma_n^*(X)$
has positive mass on the negative real axis.
%Even then, $\beta_h \leq\beta_h^*$ for all $h$.
%(iii)
%Incidentally,
%The above results are in sharp
%In contrast,
% to
%Recall that the LSD of the $S$ matrix
%whose LSD
%is supported on the
%interval $[0, 4]$. See Bai (1999).

\textup{(iii)}
Since $\Gamma_n^*(X)$ is not non-negative definite,
%the bounded Lipschitz argument of Lemma \ref{lemaux} (b) cannot be
%used. Hence
the proof of Theorem \ref{thmautocov2} for $d=\infty$ is different from
the proof of Theorem \ref{thmautocov} and needs Assumption
\ref{assumptionA}\textup{(b)}. A detailed discussion on the different
assumptions is given in Remark \ref{remautocov2} at the end of the
proofs.

\mbox{}\hspace*{1pt}\textup{(iv)} Unfortunately, the moments of the LSD of
$\Gamma_n(X)$ has no easy description. There is no easy description of the
constants $\{p_k^{(d)}\}$ either. To explain briefly
the complications involved in providing explicit expressions for these
quantities,
consider the much simpler random Toeplitz matrix $n^{-1/2}T_{n,
\varepsilon}=n^{-1/2} ((\varepsilon_{|i-j|}))$ where
$\{\varepsilon_t\}$ is i.i.d. with mean zero variance 1.
Bryc, Dembo and Jiang \cite{bry} and Hammond and Miller \cite
{hammil05} have showed that
the LSD exists and is universal. The limit moments are of the form
\[
\beta_{2k}(T)=\sum p(w),
\]
where the sum is over the so called matched words $w$ and for each $w$,
$p(w)$ is given as the volume of a suitable
subset of a $k$-dimensional hypercube. These subsets are defined
through the intersection of $k$ hyperplanes which arise from the
function $L(i,j)=|i-j|$. Thus the value of $p(w)$ can be calculated by
performing multiple integration but must be done only via numerical
integration when $k$ becomes large.
For more details, see Bose and Sen \cite{Bose08}. For our set up,
definition of matched words is generalised and is given in Section \ref
{sectionproofs} and $p_k^{(d)}$ are given by more complicated
integrals. This is the main reason why the moments of the LSD cannot be
obtained in any closed form, even when $X$ is the i.i.d. process.

Bose and Sen \cite{Bose08} considered the Toeplitz matrix $T_{n,
X}=((X_{|i-j|}))$ and showed that its LSD exists under suitable conditions.
The moments $\beta_{2k}^*$ of the LSD can be written in terms of $\{
\theta_j\}$ and $\{\beta_{2k}(T)\}$. This relation is given by
%
%e2.5 #&#
\begin{equation}
\label{eqiidvsdep} \beta_{2k}^*= \E\Biggl\llvert\sum
_{j=0}^{\infty} \theta_j \exp(-2\uppi  \RMi  j U)
\Biggr\rrvert^{2k}\beta_{2k}(T),
\end{equation}
where $U$ is uniformly distributed on $(0, 1)$.

Even a relation like (\ref{eqiidvsdep}) relating the i.i.d. process
case to the
linear process case eludes us for the autocovariance matrix.
This is primarily due to the non-linear dependence of the
autocovariances $\{\hat\gamma_X(k)\}$ on the driving
$\{\varepsilon_t\}$. One of the Referees has pointed out that in this
context, the so called ``diagram formula'' (see Arcones \cite
{Arcones}, Giraitis, Robinson and Surgailis
\cite{GRS} for details) may be useful, presumably to obtain a formula
relating the linear process case to the i.i.d. case.

It is also noteworthy that no limit moment formula or explicit
description of the LSD is known for the matrix
$n^{-1}H_{n, \varepsilon}H^\prime_{n, \varepsilon}$ where $H_{n,
\varepsilon}$ is the \textit{non-symmetric} Toeplitz matrix defined
using an i.i.d. sequence (see
Bose, Gangopadhyay and Sen \cite{bosegangosen10}).
\end{rem}

%Theorem \ref{thmautocov} hold
%mutas mutandis for $\Gamma_n^\alpha(X)$ with a suitable sequence $\{p^{

%Suppose Assumption A(a) or A(b) holds and $\alpha\ne0$.
%
% (a) Then almost surely, $F^{\Gamma_{n}^{\alpha}(X^{(d)})}
% depend on the distribution of $\{\varepsilon_{t}\}$.
% Further, for some sequence of constants $\{p^{\alpha, (d)}_{
% \begin{equation}\label{eqnautocov-limit-moment}
% (b) Under Assumption B(c), almost surely, $F^{\Gamma_{n}^{
%independent of the distribution of $\{\varepsilon_{t}\}$. Further as
%$d\to\infty$,
% $$F_d^\alpha\stackrel{w}{\to}F^\alpha\mbox{and} \beta_{h,d}^
% (c) Under Assumption B(a), $F_d^\alpha$ has unbounded support and
% $\beta_{h,d-1}^\alpha\leq\beta_{h,d}^\alpha$ if $d\geq1$.
%Hence if
% Assumption B(b) and B(c) hold, then $F^\alpha$ has unbounded
%support.
%
%(d)}_{\mathbf{k}}\}$ instead of $\{p^{(d)}_{\mathbf{k}}\}$ and under
%similar assumptions, part (a) and (b) of Theorem \ref{thmautocov}
%hold for $\Gamma_n^0(X)$. Moreover for $d$ finite or infinite the
%limit distributions have bounded support.
%

%th2.3 #&#
\begin{theorem}[(Banded and tapered sample ACVM)]\label{thmautocovband}
Suppose Assumption \textup{\ref{assumptionA}(b)} holds.

\textup{(a)} Let $0 < \alpha\leq1$. Then all the conclusions of Theorem \ref{thmautocov}
hold for $\Gamma_n^{\alpha, I}(X^{(d)})$ and $\Gamma_n^{\alpha,
\mathit{II}}(X^{(d)})$ with modified universal constants $\{p^{\alpha, I,
(d)}_{\mathbf{k}}\}$ and $\{p^{\alpha, \mathit{II}, (d)}_{\mathbf{k}}\}$,
respectively,
in (\ref{eqnautocov-limit-moment}). Same conclusions continue to hold
also for $d= \infty$.

\textup{(b)} If $\alpha=0$, and Assumption \textup{\ref{assumptionB}(b)} holds, the LSD
of $\Gamma_n^{\alpha, I}(X)$ and $\Gamma_n^{\alpha, \mathit{II}}(X)$
are $f_X(U)$.\vspace*{1pt}
%For $\alpha=0$ conclusions of part (a) and (b) of Theorem
%limit distributions
%for both $\Gamma_n^{\alpha, I}(X^{(d)})$ and $\Gamma_n^{\alpha,
%$\{F_d^{0, I}\}$, $F^{0,I}$, $\{F_d^{0, \mathit{II}}\}$ and $F^{0,\mathit{II}}$ all have
%bounded support. Moreover,
%F_d^{0, I}=F_d^{0, \mathit{II}} \mbox{ for all } d \mbox{ and }
%F^{0,I}=F^{0,\mathit{II}}.

\textup{(a)} and \textup{(b)} remain true for $\Gamma_n^{\alpha, \mathit{II}}(X^{(d)})$ and
$\Gamma_n^{\alpha, \mathit{II}}(X)$ under Assumption \textup{\ref{assumptionA}(a)}.

%Suppose Assumption A(b) and
\textup{(c)} Suppose Assumption \textup{\ref{assumptionB}(b)} holds. Let
$K$ be bounded, symmetric and continuous at 0, $K(0)=1$, $K(x)
=0$ for $|x| >1$. Suppose $m_n\rightarrow\infty$ such that $m_n/n \to
0$. Then the LSD of $\Gamma_{n, K}(X)$ is $f_X(U)$ for $d \leq\infty$.
%exists and equals is the same as that of
%$\Sigma_n(X)$ which equals
\end{theorem}

%re2.2 #&#
\begin{rem}\label{remkernel}
\textup{(i)} When $K$ is non-negative definite, Theorem
\ref{thmautocovband}\textup{(c)}
holds under Assumption \ref{assumptionA}\textup{(a)}.

\mbox{}\hphantom{i}\textup{(ii)} Xiao and Wu \cite{xiaowu} show that under the assumption
$m_n=\RMo(n^\gamma)$ (for a suitable $\gamma$)
and other conditions, the maximum eigenvalue of $\Sigma_n(X)-\Gamma
_n(X)$ tends to zero a.s.

\textup{(iii)}
%Different values of $\{\theta_k\}$ may yield
%the same LSD.
%In particular,
Each of the LSDs above
are identical for the combinations
$(\theta_0,\theta_1,\theta_2, \ldots)$,
$(\theta_0,-\theta_1,
\theta_2, \ldots)$ and $(-\theta_0,\theta_1,-\theta_2,
\ldots)$. See Basak, Bose and Sen \cite{Basakbosesen} for a
proof which is based on
properties of the limit moments.
The LSDs $f_X(U)$ of $\Sigma_n(X)$
are identical for processes with autocovariances $(\gamma_0,\gamma
_1,\ldots,\gamma_d)$ and $(\gamma_0,-\gamma_1,\ldots,(-1)^d\gamma_d)$.
The same is true of all the above LSDs.
% of $\Gamma_n(X^{(d)})$
%Same remark holds for Theorem \ref{thmautocovband}.
% and \ref{thmtaperedautocov}.
%It may be noted that
\end{rem}
\section{Proofs}\label{sectionproofs}
%Before we start the proof we wish to point out that
Szeg{\"o}'s theorem (or its triangular version)
for non-random Toeplitz matrices needs summability (or square
summability) of the entries and that is absent (in the a.s.
sense) for $\Gamma_n(X)$.
As an answer to a question raised by
Bai \cite{Bai99}, Bryc, Dembo and Jiang \cite{bry} and Hammond and
Miller \cite{hammil05} showed that
% had posed the question of the existence of the LSD for
for the random Toeplitz matrix $n^{-1/2}T_{n,
\varepsilon}=n^{-1/2} ((\varepsilon_{|i-j|}))$ where
$\{\varepsilon_t\}$ is i.i.d. with mean zero variance 1,
the LSD exists and is universal (does not depend on the
underlying distribution of $\varepsilon_1$). Bose and Sen \cite{Bose08}
considered the Toeplitz matrix $T_{n, X}=((X_{|i-j|}))$ and showed that
the LSD of $n^{-1/2}T_{n, X}$ exists under the following condition: $X$
satisfies (\ref{eqsetup}), $\sum_{j=0}^\infty|\theta_j| < \infty
$; further, $\{\varepsilon_j\}$ are independent with mean zero and
variance 1 and are (i) either uniformly bounded or (ii) are identically
distributed and $\sum_{j=0}^\infty j\theta_j^2 < \infty$.
%+is either i.i.d. with m. and with some additional assumptions.
%_t=\sum_{j=0}^\infty\theta_j \varepsilon_{t-j}$.
However, none of the above two results
%for random Toeplitz matrices
are applicable to $\Gamma_n(X)$ due to the non-linear dependence of
$\hat{\gamma}_X(k)$ on $\{X_t\}$.

Our two main tools will be (i) the moment method to show convergence of
distribution and (ii) the bounded Lipschitz metric
to reduce the unbounded case to the bounded case and also to prove the
results for the infinite order case from the
finite order case.
%Our proofs will primarily be based on
%The moment method which may be described in brief as follows.
Suppose $\{A_n\}$ is a sequence of $n \times n$ symmetric random
matrices. Let $\beta_h (A_n)$ be the $h$th moment of its ESD. It
has the following nice form:
\[
\beta_h (A_n)=\frac{1}{n} \sum
_{i=1}^n \lambda_i^h =
\frac{1}{n} \Tr\bigl(A_n^h\bigr).
\]
Then the LSD of $\{A_n\}$ exists a.s. and
%the limit distribution
is uniquely identified by its moments $\{\beta_h\}$ given below if the
following three conditions hold:

(C1) $\E[\beta_h (A_n)] \longrightarrow\beta_h$ for all $h$
(convergence of the average ESD).

(C2) $\sum_{n=1}^\infty\E[\beta_h (A_n)-\E[\beta_h
(A_n)] ]^4 < \infty$.

(C3) $\{\beta_h\}$ satisfies Carleman's condition: $\sum_{h =
1}^\infty\beta_{2h}^{-1/2h} = \infty$.

Let $d_{\mathrm{BL}}$ denote the \textit{bounded Lipschitz metric} on the space
of probability measures on $\mathbb{R}$,
%any Polish space $(X,d)$,
topologising the weak convergence of probability measures (see Dudley
\cite{Dudley}).
%$$d_{\mathrm{BL}}(\mu, \nu) = \sup\{ \int f \mrmd \mu- \int f \mrmd \nu: ||f||_{
%where $$||f||_{\infty} = \sup_x |f(x)|, ||f||_L = \sup_{x \ne y}
%|f(x)-f(y)|/d(x,y).$$
%This metric will be used to estimate the distance between spectral
%measures via the following Lemma.
%Then we have
The following lemma and its proof is given
%may be found
in
Bai \cite{Bai99}.

%le1 #&#
\begin{lem} \label{lemaux}\textup{(a)} Suppose $A$ and $B$ are $n \times n$
real symmetric matrices. Then
%
%e3.1 #&#
\begin{equation}
\label{blmetric} d^2_{\mathrm{BL}}\bigl(F^A,
F^B\bigr) %\le\left(\frac{1}{n}
\le\frac{1}{n}
\Tr(A-B)^2.
\end{equation}

\textup{(b)} Suppose $A$ and $B$ are $p \times n$ real matrices. Let $X=AA^T$
and $Y=BB^T$. Then
%
%e3.2 #&#
\begin{equation}
\label{blmetric2} d^2_{\mathrm{BL}}\bigl(F^{X},
F^{Y}\bigr) \le%\left( \frac{1}{p} \sum_{i=1}^p | \lambda_i(X) -
%^2 \le
\frac{2}{p^2} \Tr(X + Y)
\Tr\bigl[( A-B) (A-B)^T\bigr].
\end{equation}
\end{lem}
When $\alpha=1$, then without loss of generality for asymptotic
purposes, we assume that \mbox{$m_n=n$}. We visualise the full ACVM $\Gamma_n(X)$
%may be d
as the case with $\alpha=1$.
When $\{X_t\}$ is a finite order moving average process with bounded $\{
\varepsilon_t\}$, we use the \textit{method of moments} to establish
Theorem~\ref{thmautocov}\textup{(a)}. The longest and hardest part of the
proof is to verify (C1).
We first develop a manageable expression for the moments of the ESD and
then show that asymptotically only ``matched'' terms survive.
%remain in the empirical moments.
These moments are then written as an iterated sum, where one summation
is over finitely many terms (called ``words'').
Then we verify (C1) by showing that each one of these finitely many
terms has a limit.
The $d_{\mathrm{BL}}$ metric is used to remove the boundedness assumption as
well as to
deal with the infinite order case.
Easy modifications of these arguments yield the existence of the LSD
when $ 0 \leq\alpha\leq1$ in
Theorem \ref{thmautocovband}(a) and (b). The proof of Theorem \ref
{thmautocov2} is a byproduct of the arguments in the proof of Theorem
\ref{thmautocov}. However, due to the matrix now not being
non-negative definite, we impose Assumption \ref{assumptionA}(b). %the restriction that
%the random variables $\{\varepsilon_t\}$ are uniformly bounded.
The proof of Theorem \ref{thmautocov}(a) is given in details. All
other proofs are sketched and details are available in
Basak, Bose and Sen~\cite{Basakbosesen}.\looseness=-1
\subsection{Proof of Theorem \texorpdfstring{\protect\ref{thmautocov}}{2.1}}
\label{subsectionProofofTheoremautocova}

% (a)}
%($d$ finite, $m_n=n$)}
%A detailed proof will be provided only
%of Theorem \ref{thmautocov} (a) and when
%for $\alpha=1$.
%In this case,
%Without loss of generality we can take $m_n=n$.
% The proof for $\alpha\in(0,1)$ is similar and hence shall
% be only briefly sketched.
% The case $\alpha=0$ will be a bit different and therefore we will
%give a somewhat detailed proof for that situation.

The first step is to show that
we can
%may
without loss of generality, assume that $\{\varepsilon_{t}\}$ are
uniformly bounded so that we
can
%may
use the moment method.
For a standard proof of the following lemma, see
%may be found in
Basak, Bose and Sen \cite{Basakbosesen}. For convenience, we
will write
\[
\Gamma_n\bigl(X^{(d)}\bigr)=\Gamma_{n,d}.
\]

%le2 #&#
\begin{lem}\label{lembounded}
If for every
$\{\varepsilon_t\}$ satisfying Assumption \textup{\ref{assumptionA}}\textup{(b)},
$\Gamma_{n}(X^{(d)})$ has the same LSD a.s., then this
%the same
LSD continues to hold if $\{\varepsilon_t\}$ satisfies Assumption
\textup{\ref{assumptionA}}\textup{(a)}.
\end{lem}
\textit{Thus from now on we assume that Assumption} \ref{assumptionA}(b) \textit{holds.}
% Recall from Section \ref{subsectionMomentmethod} condition (C1),
%that we need to prove the convergence for every moment. Thus,
Fix any arbitrary positive integer $h$ and consider the $h$th
moment. Then
%The $h^{th}$ moment of the ESD of $\Gamma_{n,d}$ is
%
%e3.3 #&#
\begin{eqnarray}
\label{eqnbetahnd} \Gamma_{n,d}&=&\frac{1}{n}\bigl(
\bigl(Y_{i,j}^{(n)}\bigr)\bigr)_{i,j=1,\ldots, n}
\qquad\mbox{where } Y_{i,j}^{(n)}=\sum
_{t=1}^{n}X_{t,d}X_{t+|i-j|,d}\I
_{(t+|i-j| \leq n)},
\nonumber
\\
\beta_h(\Gamma_{n,d}) &=&\frac{1}{n}\Tr\bigl(
\Gamma_{n,d}^h\bigr) %\nonumber\\ &=&
= \frac{1}{n^{h+1}}\sum
_{1\leq\pi_0=\pi_h,
\pi_1,\ldots,\pi_{h-1} \leq n }Y_{\pi_0,\pi_1}^{(n)}\cdots
Y_{\pi
_{h-1},\pi_h}^{(n)}
\nonumber\\[-8pt]\\[-8pt]
&=&\frac{1}{n^{h+1}}\mathop{\sum_{1\leq\pi_0,\ldots,\pi_h \leq
n }}_{ \pi_h=\pi_0}
\Biggl[\prod_{j=1}^h \Biggl(\sum
_{t_j=1}^n X_{t_j,d}X_{t_j+|\pi
_{j-1}-\pi_{j}|,d}
\I_{(t_j+|\pi_{j-1}-\pi_{j}|\leq n)} \Biggr) \Biggr].\nonumber
\end{eqnarray}
To express the above in a neater and more amenable form, define
\begin{eqnarray*}
\mathbf{t}&=&(t_1,\ldots, t_h),\qquad \bolds{\pi}=(
\pi_0,\ldots, \pi_{h-1}),
\\
\mathcal{A}&=& \bigl\{(\mathbf{t},\bolds{\pi})\dvt  1\leq t_1,\ldots,
t_h,\pi_0,\ldots, \pi_{h-1}\leq n,
\pi_h =\pi_0 \bigr\},
\\
\mathbf{a}(\mathbf{t}, \bolds{\pi})&=& \bigl(t_1,\ldots,
t_h, t_1 + |\pi_0-\pi_1|,\ldots, t_h + |\pi_{h-1}-\pi_{h}| \bigr),
\\
\mathbf{a}& = &(a_1,\ldots, a_{2h})\in\{1,2,\ldots, 2n
\}^{2h},
\\
X_{\mathbf{a}}&=&\prod_{j=1}^{2h}(X_{a_j,d})
\quad\mbox{and}\quad \I_{\mathbf{a}(\mathbf{t}, \bolds{\pi})}=\prod_{j=1}^{h}
\I_{(t_j + |\pi_{j-1}-\pi_{j}|\leq n)}.
\end{eqnarray*}
Then using (\ref{eqnbetahnd}) we can write the so called \textit
{trace formula},
%
%e3.4 #&#
\begin{equation}
\label{eqautocov-moment}\label{eqtrace}
\E\bigl[\beta_{h}(\Gamma_{n,d})
\bigr] %&=&
% {\prod_{j=1}^{h}}X_{t_j+|\pi_j-\pi_{j-1}|,d} {
=
\frac{1}{n^{h+1}}\E\biggl[\sum_{(\mathbf{t},\bolds
{\pi}) \in\mathcal{A}}X_{\mathbf{a}(\mathbf{t}, \bolds{\pi
})}
\I_{\mathbf{a}(\mathbf{t}, \bolds{\pi})} \biggr].
\end{equation}
%
%s3.1.1 #&#
\subsubsection{Matching and negligibility of certain terms}
\label{subsubsectionmatching}
%Note that
By independence of $\{\varepsilon_t\}$,
%the expectation of the product
$\E[X_{\mathbf{a}(\mathbf{t}, \bolds{\pi})}]=0$ if there is at
least one component of the product that has no $\varepsilon_t$ common
with any other component. Motivated by this, we introduce a notion of
%appropriate
matching and
%then
show that certain higher order terms can be asymptotically neglected in
%the trace formula
(\ref{eqtrace}).
We say:

$\bullet$ $\mathbf{a}$ is \textit{$d$-matched} (in short
\emph{matched}) if $\forall i\leq2h,\exists j\neq i$ such that
$|a_i-a_j|\leq d$.
%Note that
When $d=0$ this means $a_i=a_j$.

$\bullet$ $\mathbf{a}$ is \textit{minimal $d$-matched}
(in short minimal matched) if there is a partition $\mathcal{P}$ of $\{
1,\ldots, 2h\}$,
%such that
%
%e3.5 #&#
\begin{equation}
\label{eqpartition}\{1,\ldots, 2h\}=\bigcup_{k=1}^{h}
\{i_k,j_k\},\qquad i_k<j_k
\end{equation}
such that $\{i_k\}$ are in ascending order and
\[
|a_x-a_y|\leq d \quad\Leftrightarrow\quad\{x,y\}=
\{i_k,j_k\} \qquad\mbox{for some } k.
\]
For example, for $d=1, h=3$ $(1,2,3,8, 9,10)$ is matched but not
minimal matched and $(1,2,5,6,9,10)$ is both matched and minimal matched.
%
%le3 #&#
\begin{lem}\label{lemorder-reduction}
$\#\{\mathbf{a}\dvt
\mathbf{a}$ is matched but not minimal matched$\} =\RMO(n^{h-1})$.
\end{lem}
\begin{pf}
Consider the graph with
% $2h$
vertices $\{1, 2,\ldots, 2h\}$.
%and
Vertices $i$ and $j$ have an
%with an
edge if $|a_i-a_j|\leq d$. Let $k=\#$ connected components. Consider a
typical $\mathbf{a}$.
%matched but not minimal matched.
Let $l_j$ be
the number of vertices in the $j$th component.
%= \# \mbox{vertices in the} j\mbox{th component}.$$
Since $\mathbf{a}$ is matched, $l_j\geq2$ for all $j$ and $l_{j} >2$
for at least one $j$.
Hence, $2h=\sum_{j=1}^k l_j>2k$. That implies
% \mbox{which implies}
$ k\leq h-1$.
Also if $i$ and $j$ are in the same connected component then
$|a_i-a_j|\leq2dh$. Hence, the number of $a_i$'s such that $i$ belongs
to any given component is
$\RMO(n)$ and the result follows.
\end{pf}
Now we can rewrite (\ref{eqautocov-moment}) as
\begin{eqnarray*}
\E\bigl[\beta_{h}(\Gamma_{n,d})\bigr] &=&\frac{1}{n^{h+1}}\E
\biggl[\sum_{1}X_{\mathbf{a}(\mathbf{t},\bolds
{\pi})}\I_{\mathbf{a}(\mathbf{t},\bolds{\pi})}
\biggr] +\frac{1}{n^{h+1}}\E\biggl[\sum_{2}X_{\mathbf{a}(\mathbf
{t},\bolds{\pi
})}
\I_{\mathbf{a}(\mathbf{t},\bolds{\pi})}\biggr]
\\
&&{} +\frac{1}{n^{h+1}}\E\biggl[\sum_{3}X_{\mathbf{a}(\mathbf
{t},\bolds{\pi})}
\I_{\mathbf{a}(\mathbf{t},\bolds{\pi})}\biggr]=T_1 + T_2 + T_3
\qquad\mbox{(say)},
\end{eqnarray*}
where
% $ {\sum_i, i=1, 2, 3}$ are
the three summations are
%taken
over
%all
$(\mathbf{t},\bolds{\pi})\in\mathcal{A}$ such that
$\mathbf{a}(\mathbf{t},\bolds{\pi})$ is, respectively, (i) minimal matched,
(ii) matched but not minimal matched and (iii) not matched.

By mean zero assumption, $T_3=0$. Since $X_i$'s are uniformly bounded, by
Lemma \ref{lemorder-reduction},
$T_2 \leq\frac{C}{n}$ for some constant $C$. So provided the limit exists,
%
%e3.6 #&#
\begin{equation}
\label{eqmomexp1} \lim_{n \rightarrow\infty} \E\bigl[
\beta_{h}(\Gamma_{n,d})\bigr] =\lim_{n \rightarrow\infty}
\frac{1}{n^{h+1}} \E\biggl[\mathop{\sum_{
(\mathbf{t},\bolds{\pi})\in\mathcal{A}\dvt
\mathbf{a}(\mathbf{t},\bolds{\pi})\ \mathrm{is}}}_{\mathrm{minimal}
\ \mathrm{matched}}X_{\mathbf{a}(\mathbf{t},\bolds{\pi})}
\I_{\mathbf
{a}(\mathbf{t},\bolds{\pi})}\biggr].
\end{equation}
\textit{Hence, from now our focus will be only on minimal matched words.}

%s3.1.2 #&#
\subsubsection{Verification of \textup{(C1)} for
Theorem \texorpdfstring{\protect\ref{thmautocov}}{2.1}\textup{(a)}}
\label{subsubsectionVerificationofC1}% ($d$
%finite, $m_n=n$)
%As mentioned earlier,
This is the hardest and lengthiest part of the proof. One can give a
separate and easier proof for the case $d=0$. However, the proof for
general $d$ and for
%the simpler case of
$d=0$ are developed in parallel since this helps to relate the limits
in the two cases.
%for general $d$ to limit for $d=0$.

%The idea behind the proof of (C1) is as follows.
Our starting point is equation (\ref{eqmomexp1}). We first define an
\textit{equivalence relation}
on the set of minimal matched $\mathbf{a}=\mathbf{a}(\mathbf{t},
\bolds{\pi})$. This
yields
finitely many equivalence classes.
Then we
%may
can write the sum
in (\ref{eqmomexp1}) as an iterated sum where the outer sum is over
the equivalence classes.
%(see (\ref{eqiteratedsum})).
Then we show that for every fixed equivalence class, the inner sum has
a limit.

To define the equivalence relation, consider the collection of
$(2d+1)h$ symbols (letters)
\[
\mathcal{W}_h=\bigl\{w_{-d}^k,\ldots,
w_0^k,\ldots, w_d^k\dvt  k=1,\ldots, h \bigr\}.
\]
Any minimal $d$ matched
$\mathbf{a}=(a_1,\ldots, a_{2h})$
induces a
partition as given in (\ref{eqpartition}).
%$$\mathcal{P}=\bigcup_{k=1}^{h} \{i_k,j_k \} \mbox{with} i_k < j_k
% \mbox{of} \{1,\ldots, 2h\}$$ where the $\{i_k\}$ are arranged
%in \textit{ascending order}.
With this $\mathbf{a}$, associate the \textit{word} $w=w[1]w[2]
\cdots w[2h]$ of length $2h$ where
%
%e3.7 #&#
\begin{equation}
\label{eqworddefine} w[i_k]=w_0^k,\qquad
w[j_k]=w_l^k \qquad\mbox{if }
a_{i_k}-a_{j_k}=l, 1\leq k \leq h.
\end{equation}
As an example, consider $d=1, h=3$ and
$\mathbf{a}=(a_1,\ldots, a_6)=(1,21,1,20,39,40)$. Then the unique
partition of $\{1, 2,\ldots, 6\}$ and the unique word associated with
$\mathbf{a}$ are
$\{\{1, 3\}, \{2, 4\}, \{5, 6\}\}$ and
$[w_0^1w_0^2w_0^1w_{1}^2w_0^3w_{-1}^3]$, respectively.

%It is important to
Note that corresponding to any fixed partition ${\cal P}=\{\{i_k,j_k\},
1\leq k \leq h\}$, there are several $\mathbf{a}$ associated with it
and there are exactly $(2d+1)^h$ words that can arise from it. For
example, with $d=1, h =2$ consider the partition ${\cal P}= \{\{1,2\},\{
3,4\}\}$. Then the nine words corresponding to ${\cal P}$ are
$w_0^1w_i^1w_0^2w_j^2$ where $i, j=-1, 0, 1$.

By a slight abuse of notation, we write $w \in{\cal P}$ if the
partition corresponding to $w$ is same as~${\cal P}$.
We will say that:

$\bullet$
$w[x]$ \textit{matches} with $w[y]$ (say
%denote by
$w[x] \approx w[y]$) iff $w[x]=w^k_l$ and $w[y]=w^k_{l'}$ for some
$k,l,l'$.

$\bullet$
$w$ is $d$ \textit{pair matched} if it is \textit{induced} by a
minimal $d$ matched $\mathbf{a}$ (so
$w[x]$ matches with $w[y]$ iff $|a_x-a_y|\leq d$).

This induces an \textit{equivalence relation} on all $d$ minimal
matched $\mathbf{a}$ and the equivalence classes can be indexed by $d$
pair matched $w$.
Given such a
%a $d$ pair matched
$w$, the corresponding equivalence class is given by
%
%e3.8 #&#
\begin{eqnarray}
\label{eqeqclass}\label{eqpiw}
\Pi(w)&=& \bigl\{(\mathbf{t},\bolds{\pi})\in\mathcal{A}\dvt
w[i_k]=w_0^k, w[j_k]=w_l^k
\nonumber\\[-8pt]\\[-8pt]
&&\hspace*{5.7pt} \Leftrightarrow\mathbf{a}(\mathbf{t}, \bolds{\pi})_{i_k}-a(\mathbf{t},
\bolds{\pi})_{j_k}=l \mbox{ and } \I_{\mathbf{a}(\mathbf
{t},\bolds{\pi})}=1\bigr\}.\nonumber
\end{eqnarray}
Then we
%may
rewrite (\ref{eqmomexp1}) as (provided the second limit exists)
%
%e3.9 #&#
\begin{equation}
\label{eqiteratedsum}\lim_{n \rightarrow\infty} \E\bigl[
\beta_{h}(\Gamma_{n,d})\bigr] = \sum
_{\cal P}\sum_{{w \in\cal P}} \lim
_{n \rightarrow\infty
}\frac{1}{n^{h+1}}\sum_{(t, \pi) \in\Pi(w)}
\E[X_{\mathbf
{a}(\mathbf{t},\bolds{\pi})}\I_{\mathbf{a}(\mathbf{t},\bolds
{\pi})}].
\end{equation}
By using the autocovariance structure, we further simplify the above as
follows. Let
\[
{\cal W}(\mathbf{k})=\bigl\{w\dvt  \#\bigl\{s\dvt  \bigl|w[i_s]-w[j_s]\bigr|=i
\bigr\}=k_i, i=0, 1,\ldots, d\bigr\}.
\]
Using the definitions of $\gamma_{X^{(d)}}(\cdot)$ and of $S_{h,d}$
given in (\ref{eqkShd}), we
%may
rewrite (\ref{eqiteratedsum}) as (for any set $Z$, $\# Z$ denotes the
number of elements in $Z$)
%
%e3.10 #&#
\begin{equation}
\label{eqmom} \lim_{n \rightarrow\infty} \E\bigl[\beta_{h}(
\Gamma_{n,d})\bigr] = \sum_{{\cal P}} \sum
_{S_{h,d}}\sum_{w \in{\cal P}
\cap{\cal W}(\mathbf{k})}\lim
_{n \rightarrow\infty} \frac
{1}{n^{h+1}}\#\Pi(w) \prod
_{i=0}^d\bigl[\gamma_{X^{(d)}}(i)
\bigr]^{k_i}
\end{equation}
provided the following limit exists for every word $w$ of length $2h$.
%
%e3.11 #&#
\begin{equation}
\label{eqdefinepw} p_w^{(d)} \equiv\lim
_{n \rightarrow\infty} \frac{1}{n^{h+1}} \#\Pi(w).
\end{equation}
To show that this limit exists, it is convenient to work with
$\Pi^{*}(w) \supseteq\Pi(w)$ defined as
%
%e3.12 #&#
\begin{eqnarray}\label{eqpi*w}
\Pi^\ast(w)&=&\bigl\{(\mathbf{t}, \bolds{\pi})\in\mathcal{A}\dvt
w[i_k]=w_0^k, w[j_k]=w_l^k
\nonumber\\[-8pt]\\[-8pt]
&&\hspace*{5.7pt}\Rightarrow a(\mathbf{t}, \bolds{\pi})_{i_k}-a(\mathbf{t}, \bolds{
\pi})_{j_k}=l \mbox{ and } \I_{\mathbf{a}(\mathbf{t},\bolds{\pi
})}=1\bigr\}.
\nonumber
\end{eqnarray}
By Lemma \ref{lemorder-reduction}, we have for every $w$,
${n^{-(h+1)}}\#(\Pi^\ast(w)-\Pi(w))\to0$.
Thus, it is enough to show that $\lim_{n \rightarrow\infty} \frac
{1}{n^{h+1}} \#\Pi^{*}(w)$ exists.

For a pair matched $w$, we divide its coordinates according to the
position of the matches as follows. For $1 \le i < j \le h$, let the
sets $S_i$ be defined as
\begin{eqnarray*}
S_1(w)&=&\bigl\{i\dvt  w[i] \approx w[j]\bigr\},\qquad
S_2(w)= \bigl\{ j\dvt  w[i] \approx w[j] \bigr\},
\\
S_3(w)&=& \bigl\{i\dvt  w[i] \approx w[j+h]\bigr\},\qquad
S_4(w)= \bigl\{j\dvt  w[i] \approx w[j+h]\bigr\},
\\
S_5(w)&=& \bigl\{i\dvt  w[i+h] \approx w[j+h]\bigr\},\qquad
S_6(w)= \bigl\{j\dvt  w[i+h] \approx w[j+h]\bigr\}.
\end{eqnarray*}
Let $E$ and $G\subset E$ be defined as
%$\begin{eqnarray*}E&=&\{t_1,\ldots, t_h, \pi_0,\ldots, \pi_h\}, \\
%G&=&\{t_i| i \in S_1(w) \cup S_3(w)\} \cup\{\pi_0\} \cup\{\pi_i | i+h
%
\begin{eqnarray*}
E&=&\{t_1,\ldots, t_h, \pi_0,\ldots,
\pi_h\},\\
G&=&\bigl\{t_i| i \in S_1(w) \cup
S_3(w)\bigr\} \cup\{\pi_0\} \cup\bigl\{
\pi_i | i+h \in S_5(w)\bigr\}.%\subset E.
\end{eqnarray*}
%
%Note that $G\subset E$.
Elements in $G$ are the indices where any matched letter appears for
the \textit{first} time
and these will be called the \textit{generating vertices}.
%These are the indices where the first occurrence of any matched letter
%happens.
%Note that
$G$ has $(h+1)$ elements say $u_1^n,\ldots, u_{h+1}^n$ and for
simplicity we will write
\[
G\equiv U_n=\bigl(u_1^n,\ldots,
u_{h+1}^n\bigr) \quad\mbox{and}\quad {\cal N}_n=\{
1, 2,\ldots, n\}.
\]

\begin{claim}\label{claim1}
Each element of $E$
%may be written as
is a  linear expression (say $\bolds{\lambda}_i$) of the generating
vertices that are all to the \emph{left} of the element.
\end{claim}
\begin{pf} Let the constants in the proposed linear expressions be $\{
m_j\}$.\vadjust{\goodbreak} %and %the linear combinations by $\{\bolds{\lambda}_i\}$.

(a) For those elements of $E$ that are generating vertices, we take the
constants as $m_j=0$ and the linear combination is taken as the
identity mapping so that
\begin{eqnarray*}
\mbox{for all } i \in S_1(w) \cup S_3(w)\qquad \bolds{
\lambda}_i&\equiv& t_i,
\\
\bolds{\lambda}_{h+1} &\equiv& \pi_0,
\end{eqnarray*}
and for all
\[
i+h \in S_5(w),\qquad \bolds{\lambda}_{i+h+1} \equiv \pi_i.
\]

(b) Using the relations between $S_1(w)$ and
$S_2(w)$ induced by $w$, we can write
\[
\mbox{for all } j\in S_2(w)\qquad t_j=\bolds{
\lambda}_j+n_j
\]
for some $n_j$ such that $|n_j| \le d$ and define $m_j=n_j$ for $j \in
S_2(w)$ and $\bolds{\lambda}_j\equiv\bolds{\lambda}_i$.

(c) Note that
for every $\bolds{\pi}$ we can write
\[
|\pi_{i-1}-\pi_i|=b_i(\pi_{i-1}-
\pi_i) \qquad\mbox{for some } b_i \in\{-1,1\}.
\]
Consider the vector
$\mathbf{b}=(b_1, b_2,\ldots, b_h) \in\{-1,1\}^h$. It will be a
valid choice if we have
%
%e3.13 #&#
\begin{equation}
\label{eqnnonneg1}b_i(\pi_{i-1}-
\pi_i) \geq0 \qquad\mbox{for all } i.
\end{equation}
We then have the following two cases:

Case 1: $w[i]$ matches with $w[j+h]$, $j+h \in S_4(w)$ and $i \in
S_3(w)$. Then we get
%one of the following two:
%
%e3.14 #&#
\begin{equation}
\label{eqnj1} t_i=t_j+b_j(
\pi_{j-1}-\pi_j)+n_{j+h} \qquad\mbox{for some
integer } n_{j+h} \in\{-d,\ldots, 0,\ldots, d\}.
\end{equation}

Case 2:
%Note that if
$w[i+h]$ matches with $w[j+h]$,
% then for
$j+h \in S_6(w)$ and $i+h \in S_5(w)$. Then we have
%
%e3.15 #&#
\begin{equation}
\label{eqnj2} t_i + |\pi_{i-1} - \pi_i| =
t_j + |\pi_{j-1} - \pi_j|+n_{j+h}
\qquad\mbox{where } n_{j+h} \in\{-d,\ldots, 0,\ldots, d\}.
\end{equation}
So we note that inductively from left to right we can write
%
%e3.16 #&#
\begin{equation}
\label{eqmj}\pi_j=\bolds{\lambda}^{\mathbf{b}}_{j+1+h}+m_{j+1+h},\qquad
j+h \in S_4(w) \cup S_6(w).
\end{equation}
Hence, inductively, $\pi_j$ as a linear combination $\{\bolds
{\lambda}^{\mathbf{b}}_{j}\}$ of the generating vertices up to
an appropriate constant.
The superscript $\mathbf{b}$
emphasizes
that $\{\bolds{\lambda}^{\mathbf{b}}_{j}\}$
%the linear expression will
depends on $\mathbf{b}$. Further,
%Also by construction,
$\{\bolds{\lambda}^{\mathbf{b}}_{j}\}$
%for any element
depends \textit{only} on the vertices present to the left of it. %This
%proves our claim.
\end{pf}

Now we are almost ready to write down an expression for the limit. If
$\bolds{\lambda}_i$ were unique for each~$\mathbf{b}$, then we
could write $\#\Pi^*(w)$ as a sum of all possible choices of $\mathbf
{b}$ and we could tackle the expression for each $\mathbf{b}$
separately. However, $\bolds{\lambda}_i$'s may be same for several
choices $b_i \in\{-1,1\}$. For example, for the word
$w_0^1w_0^2w_0^1w_0^2$, we can
%may
choose any $\mathbf{b}$.\vadjust{\goodbreak} We
%first
circumvent this problem as follows:
%in the following way:
%So we need to be careful about using ${\mathbf b}$ as an index of
%summation. To tackle this,
Let
\[
{\cal T}= \bigl\{ j+h \in S_4(w) \cup S_6(w)|
\bolds{\lambda}^{\mathbf{b}}_{j+h} - \bolds{\lambda}^{\mathbf{b}}_{j+h-1}
\equiv0\ \forall b_j \bigr\}.
\]
Note that the definition of ${\cal T}$ depends on $w$ only through the
partition ${\cal P}$ it generates.

Suppose $j+h \in{\cal T}$.
Define
%
%e3.17 #&#
%e3.18 #&#
\begin{eqnarray}
L_j(U_n) &:=& b_j\bigl(\bolds{\lambda}
^{\mathbf{b}}_{j+h-1}\bigl(U^n\bigr)-\bolds{\lambda
}^{\mathbf{b}}_{j+h}\bigl(U^n\bigr)\bigr)+m_{j+h-1}-m_{j+h}
\\
&:= & \tilde L_j(U_n)+m_{j+h-1}-m_{j+h}.
\end{eqnarray}
Then from (\ref{eqnj1}) and (\ref{eqnj2}) the region given by (\ref
{eqnnonneg1}) is
%
%e3.19 #&#
\begin{equation}
\bigl\{L_j(U_n) \geq0\bigr\}\equiv\bigl\{ \tilde
L_j(U_n)+m_{j+h-1}-m_{j+h}\geq0 \bigr
\}.
\end{equation}

\begin{claim}\label{claim2}
The above expression is same for all choices of $\{
b_j\}$, for $j + h \in{\cal T}$.
\end{claim}

%Here is a short proof of the claim.
%
\begin{pf} First, we show that if $j+h \in{\cal T}$ then we must have
%
%e3.20 #&#
\begin{equation}
\label{eqtau} t_j = t_j +|\pi_{j-1} -
\pi_j|+n_j \qquad\mbox{for some integer } |n_j|
\leq d.
\end{equation}
Suppose this is not true.
So first assume that $j+h \in S_6(w)$. Then we will have a relation
%
%e3.21 #&#
\begin{equation}
\label{eqclaim2rel}t_i + b_i (\pi_{i-1} -
\pi_i) = t_j + b_j(\pi_{j-1} -
\pi_j) + n_j\qquad \mbox{where } i+h \in
S_5(w).
\end{equation}
%
%Recall that any typical linear function,
Since $\bolds{\lambda}^{\mathbf{b}}_{j}$ depends only on the vertices present
to the left of it,
%. Thus
in (\ref{eqclaim2rel}),
%the above equation
coefficient of $\pi_i$ would be non-zero and hence we must have
$\bolds{\lambda}_{j+h-1}^{\mathbf{b}} - \bolds{\lambda}_{j+h}^{\mathbf{b}} \not
\equiv0$.\vspace*{1pt}

Now assume $j+h \in S_4(w)$ and $w[i]$ matches with $w[j+h]$ for $i \ne
j$. Then we can repeat the argument above to arrive at a similar
contradiction. This shows that if $j+h \in{\cal T}$ then our relation
must be like (\ref{eqtau}). Now a simple calculation shows that for
such relations,
%like (\ref{eqtau}) we have
%
\[
b_j\bigl(\bolds{\lambda}^{\mathbf{b}}_{j+h-1}(U_n)-
\bolds{\lambda}^{\mathbf{b}}_{j+h}(U_n)
\bigr)+m_{j+h-1}-m_{j+h}= -n_j,
\]
which is of course same across all choices of $\mathbf{b}$. This
proves our claim.
\end{pf}

Now note that if $j+h \in{\cal T}$ and if $n_{j+h} \neq0$ then as we
change $b_j$ it does change the value of $m_{2h+1}$. Further, we can
have at most two choices for $\pi_{j}$ for every choices of $\pi
_{j-1}$ if $n_{j+h} \ne0$ depending on $b_j$.

However for $j+h \in{\cal T}$ and $n_j =0$,
%none of them is true and
we have only one choice for $\pi_j$ given the choice for $\pi_{j-1}$
for every choice of $b_j$. On the other hand, we know $\mathbf{b} \in
\{-1,1\}^h$ must satisfy (\ref{eqnnonneg1}).
Keeping the above in view,
let
\[
{\cal B}(w) = \bigl\{ \mathbf{b} \in\{-1,1\}^h|
b_j=1 \mbox{ if } n_j=0 \mbox{ for } j \in{
\cal T} \bigr\},\vadjust{\goodbreak}
\]
where $\{n_j\}$ is as in Claim \ref{claim2}. For ease of writing, we introduce a
few more notation:
%
%e3.22 #&#
\begin{eqnarray}\label{eqindicators}
\I_{m,h}(U_n) &:=& \I\bigl(\bolds{\lambda}^{\mathbf
{b}}_{2h+1}(U_n)+m_{2h+1}
= \bolds{\lambda}^{\mathbf
{b}}_{h+1}(U_n)+m_{h+1}
\bigr),
\nonumber
\\
\I_{\bolds{\lambda}^{\mathbf{b}}, L} (U_n)&:= & \prod_{j=1}^{h}
\I\bigl(\bolds{\lambda}^{\mathbf{b}}_{j}(U_n)+L_j(U_n)
\leq n\bigr),\nonumber\\[-8pt]\\[-8pt]
\I_{\bolds{\lambda}^{\mathbf{b}},m}(U_n)&:=&\prod
_{j=1}^{2h}\I\bigl(\bolds{\lambda}
^{\mathbf{b}}_j(U_n)+m_j\in{\cal
N}_n\bigr)\quad\mbox{and}
\nonumber\\
\I_{{\cal T}} (U_n)&:= &\prod
_{1\leq j \leq h, j
\notin{\cal T}}\I\bigl(L_j(U_n) \geq0\bigr)
\times\prod_{j \in{\cal T}} \I( n_j \leq0).
\nonumber
\end{eqnarray}
Now we note that,
\begin{eqnarray*}
p_w^{(d)} &:= &\lim_n
\frac{1}{n^{h+1}} \#\Pi^\ast(w)
\\
& = &\lim_n\frac{1}{n^{h+1}} \sum
_{\mathbf{b}\in{\cal B}(w)}\sum_{U_n\in{\cal N}_n^{h+1}} \I
_{m,h}(U_n) \times\I_{\bolds{\lambda}^{\mathbf{b}},m}(U_n)
\times\I_{\bolds{\lambda}^{\mathbf{b}},L}(U_n)\times\I_{{\cal T}}
(U_n)
\\
& =& \lim_n \sum_{\mathbf{b}\in{\cal B}(w)}
\E_{U_n} \bigl[ \I_{m,h}(U_n) \times
\I_{\bolds{\lambda}^{\mathbf{b}},m}(U_n) \times\I_{\bolds{\lambda
}^{\mathbf{b}},L}(U_n)
\times\I_{{\cal T}} (U_n) \bigr].
\end{eqnarray*}
Now it only remains to identify the limit. To this end, first fix a
partition ${\cal P}$ and $\mathbf{b} \in\{-1,1\}^h$. If $d=0$, then
there is one and only one word corresponding to it.
However, across any $d$ and any fixed $k_0, k_1,\ldots, k_d$, the
linear functions $\bolds{\lambda}_j$'s continue to remain same. The
only possible changes will be in the values of $m_j$'s. %\textit{This
%is why there is a relation between the limit for $d=0$ and $d\ne0$}.

We now identify the cases where the above limit is zero.

\begin{claim}\label{claim3}
Suppose $w$ is such that ${\cal R}: =
\{\lambda^{\mathbf{b}}_{2h+1}(U_n)+m_{2h+1} =\bolds{\lambda}
^{\mathbf{b}}_{h+1}(U_n)+m_{h+1} \}$ is a lower dimensional
subset of ${\cal N}_n^{h+1}$. Then
the above limit is zero.
\end{claim}

\begin{pf}
First, consider the case $d=0$. Then $m_j=0,\forall j$.
%in that case, in the expressions above the $m_j$'s are absent. Fix a
%word $w$ (or equivalently a partition ${\cal P}$).
Note that
${\cal R}$ lies in a hypercube. Hence, the result follows by
convergence of
the Riemann sum to the corresponding Riemann integral.
For any general $d$, the corresponding region is just a translate of
the region considered for $m_j=0$.
Hence, the result follows.
\end{pf}

Hence for a fixed $w \in{\cal P}$, a positive limit contribution is possible
only when ${\cal R}= {\cal N}_n^{h+1}$. This implies that we must have
\begin{eqnarray*}
\bolds{\lambda}^{\mathbf{b}}_{2h+1}(U_n) -\bolds{
\lambda}^{\mathbf{b}}_{h+1}(U_n) &\equiv&0 \qquad\mbox{(for
$d=0$)},
\\
\bolds{\lambda}^{\mathbf{b}}_{2h+1}(U_n) -\bolds{\lambda}
^{\mathbf{b}}_{h+1}(U_n) &\equiv&0 \quad\mbox{and}\quad
m_{2h+1}-m_{h+1}=0 \qquad\mbox{(for general $d$)}.
\end{eqnarray*}
Note that the first relation depends only the partition ${\cal P}$ but
the second relation
is determined by the word $w$.
%Now fix a partition ${\cal P}$ such that $\bolds{\lambda^{
%0$ then
%$$p_{n,w}^{(d)}:= \frac{1}{n^{h+1}} |\Pi^{*}(w)|=0.$$ Thus the limit
%contribution $p_w^d$ from this word will be again zero. Using these
%facts
%We are now ready to write an expression for the limit.
%???I have suppressed a few lines here Ok to do that??
Now $\bolds{\lambda}^{\mathbf{b}}_j$ being linear forms with integer coefficients
\[
\bolds{\lambda}^{\mathbf{b}}_j(U_n)+m_j
\in\{1,\ldots, n\}\quad\iff\quad\bolds{\lambda}^{\mathbf{b}}_j \biggl(
\frac{U_n}{n} \biggr)+\frac{m_j}{n}\in(0,1].
\]
Define $\I_{m,h}(U)$, $\I_{\bolds{\lambda}^{\mathbf{b}}, \tilde
{L}} (U)$, $\I_{\bolds{\lambda}^{\mathbf{b}}}(U)$
and $\tilde{\I}_{{\cal T}} (U)$ as in
(\ref{eqindicators}) with $U_n$ replaced by $U$, $L$ replaced by
$\tilde L$, ${\cal N}_n$ replaced by $(0,1)$, $n$ replaced by $1$, and
dropping $m_j$'s in $\I_{\bolds{\lambda}^{\mathbf{b}},m}$.
% \I_{m,h}(U)&:=& \I(\bolds{\lambda^{\mathbf{b}}}_{2h+1}(U)= \bolds{
% \mbox{ and } & & \tilde{\I}_{{\cal T}} (U):= \prod_{{{j=1} \atop{j
Noting $\frac{U_n}{n}\stackrel{w}{\Rightarrow}U$ following uniform
distribution on $[0,1]^{h+1}$, $\frac{1}{n^{h+1}}\lim{}\#\Pi^\ast(w)$
%=p_w^{(d)}$
equals
%
%e3.23 #&#
\begin{equation}
p_w^{(d)} = \sum_{ \mathbf{b} \in{\cal B}(w)}
\E_U \bigl[ \I_{m,h}(U) \times\I_{\bolds{\lambda}^{\mathbf{b}}, \tilde
{L}} (U)
\times\I_{\bolds{\lambda}^{\mathbf{b}}} (U) \times\tilde{\I}_{{\cal
T}}(U) \bigr].
\end{equation}
Now the verification of (C1) is complete by observing that (\ref
{eqmom}) becomes
\begin{eqnarray}
\label{eqmom2}
\lim_{n \rightarrow\infty} \E\bigl[\beta_{h}(
\Gamma_{n,d})\bigr] &=& \sum_{\cal P} \sum
_{\mathbf{k} \in S_{h,d}}p_{\mathbf
{k}}^{{\cal P}, d} \prod
_{i=0}^d\bigl[\gamma_{X^{(d)}}(i)
\bigr]^{k_{i}} \nonumber\\[-8pt]\\[-8pt]
&=& \sum_{\mathbf{k} \in S_{h, d}}
p_{\mathbf{k}}^{(d)}\prod_{i=0}^d
\bigl[\gamma_{X^{(d)}}(i)\bigr]^{k_i},\nonumber
\end{eqnarray}
where
%
%e3.24 #&#
\begin{equation}
\label{eqlimitpw} p_{\mathbf{k}}^{{\cal
P},d}=\sum
_{ w \in{\cal P} \cap{\cal W}(\mathbf{k})} p_w^{(d)}
\quad\mbox{and}\quad
p^{(d)}_{ \mathbf{k}}= \sum_{{\cal P}}
p_{\mathbf{k}}^{{\cal P},d}.
\end{equation}
%
%From the above discussion, observing that the indicators in (
% it follows that
%$$p_w^{(d)} \leq p_w^{0} 2 ^h.$$
%As a consequence we have
%p_{\mathbf{k}}^{\mathcal P, d} \leq2^h \Bigg({{h}\atop{k_0, k_1,
%
Since there is no explicit expression for the moments of the LSD, we
provide in Table~\ref{tablemoments} the first three moments of the
LSD of $\Gamma_n(X)$, when the input sequence is i.i.d. and MA(1). To
calculate the moments, we need to find the contributions $p_w^{(d)}$
for words $w$. The contributions of different relevant words, are
provided in Table \ref{tablemoment2}, and in Table \ref
{tablemoment3}, for the i.i.d. case. For the MA(1), one can work out
the contributions from there.

%t1 #&#
\begin{table}
\tablewidth=150pt
\caption{Contributions from words of length $4$ for i.i.d.
case}\label{tablemoment2}\vspace*{-1pt}
\begin{tabular*}{\tablewidth}{@{\extracolsep{\fill}}ll@{}}
\hline
Word $w$ & Contribution $p_w^{(0)}$ \\
\hline
aabb & $2/3$ \\
abab & 1 \\
abba & 0 \\
\hline
\end{tabular*}  \vspace*{-5pt}
\end{table}
%

%Partition $w$ & Contribution $p_w^{(0)}$ \\
%aabccb& 2/3 \\
%aabbcc & 1/6 \\
%aabcbc & 1/6 \\
%ababcc & 1/6\\
%abacbc & 2/3\\
%abaccb& 1/6\\
%abbacc & 2/3\\
%abbcca & 1/6\\
%abbcac & 1/6\\
%abcabc & 1\\
%abcacb & 0 \\
%abcbac & 0 \\
%abcbca & 0 \\
%abccab & 0 \\
%abccba & 0 \\
%
%t2 #&#
\begin{table}
\tablewidth=270pt
\caption{Contributions from words of length $6$ for i.i.d. case}\label{tablemoment3}\vspace*{-1pt}
\begin{tabular*}{\tablewidth}{@{\extracolsep{\fill}}llll@{}}
\hline
Word $w$ & Contribution $p_w^{(0)}$ & Word $w$ & Contribution
$p_w^{(0)}$ \\
\hline
aabccb& $2/3$ &abbcac &$1/6$ \\
aabbcc & $1/6$ & abcabc& 1\\
aabcbc & $1/6$ & abcacb& 0 \\
ababcc & $1/6$ & abcbac &0 \\
abacbc & $2/3$ & abcbca &0 \\
abaccb& $1/6$ & abccab&0 \\
abbacc & $2/3$ & abccba&0 \\
abbcca & $1/6$ & & \\
\hline
\end{tabular*}   \vspace*{-5pt}%
\end{table}
%

%t3 #&#
\begin{table}
\tablewidth=270pt
\caption{First three moments for i.i.d. and MA(1) input sequence}\label{tablemoments}\vspace*{-1pt}
\begin{tabular*}{\tablewidth}{@{\extracolsep{\fill}}lll@{}}
\hline
& i.i.d. & MA(1)\\
\hline
Mean & $\theta_0^2$ & $\theta_0^2+\theta_1^2$ \\[3pt]
Second moment & $\frac{5}{3} \theta_0^4$ & $\frac{5}{3}(\theta
_0^2+\theta_1^2)^2 +\frac{20}{3}\theta_0^2 \theta_1^2$ \\[3pt]
Third moment & $ 4\theta_0^6$ & $4 (\theta_0^2+\theta_1^2)^3 + 24
(\theta_0^2+\theta_1^2) (2 \theta_0 \theta_1)^2 $ \\
\hline
\end{tabular*}    \vspace*{-3pt}
\end{table}
%

%s3.1.3 #&#
\subsubsection{Verification of \textup{(C2)} and \textup{(C3)}
for Theorem \texorpdfstring{\protect\ref{thmautocov}}{2.1}\textup{(a)}}
\label{subsubsectionVerificationofC2} %($d$ finite, $m_n=n$)
%We state this formally in the following Lemma.
%The following Lemma implies (C2) and (C3).
%
%le4 #&#
\begin{lem}\label{lemalmostsure}
\textup{(a)} $\E[n^{-1}\Tr(\Gamma_{n,d}^h)-n^{-1}\E
[\Tr(\Gamma_{n,d}^h)] ]^4=\RMO (n^{-2} )$.
Hence $\frac{1}{n}\Tr(\Gamma_{n,d}^h)$ converges to $\beta_{h,d}$
a.s.

\textup{(b)}
%The sequence
$\{\beta_{h,d}\}_{h\geq0}$ satisfies \textup{(C3)} and hence defines a
unique probability distribution on $\mathbb{R}$.
\end{lem}
\begin{pf}
Proof of part (a) uses ideas from Bryc, Dembo and Jiang \cite{bry} but
%one needs to argue slightly differently as
the inputs of the matrix are no longer independent, and therefore some
modifications are needed.
%We omit the
Details are available in Basak, Bose and Sen \cite
{Basakbosesen}.

(b) Using\vspace*{2pt}
(\ref{eqmom2}) and (\ref{eqnautocov-limit-moment})
%for $p^{(d)}_{\mathbf{k}}$ and $\beta_{h,d}$ and from
%Now the
and noting that the number of ways of choosing the partition
$\{1,\ldots,2h\}=\bigcup_{l=1}^{h} \{i_l,j_l \}$ for
$\mathbf{a}(\mathbf{t},\bolds{\pi})$ is $\frac{(2h)!}{2^h h!}$,
it easily follows that
%.
%
%We recall the formulae for
% $p^{(d)}_{\mathbf{k}}$ and $\beta_{h,d}$ from (\ref{eqmom2})
%and (\ref{eqnautocov-limit-moment}).
% Now the number of ways of choosing the partition $ {
% $\mathbf{a}(\mathbf{t},\bolds{\pi})$ is $\frac{(2h)!}{2^h h!}$.
%Hence
%$$p^{(d)}_{\mathbf{k}}\leq{\lim_n}\frac{1}{n^{h+1}}
%Hence we have,
%
%e3.25 #&#
\begin{eqnarray}\label{eqsupbeta}
|\beta_{h,d}|&\leq& \sum_{S_{h,d}}
\frac{4^h(2h)!}{h!}\frac
{h!}{k_0! \cdots k_d!}\prod_{i=0}^d
\bigl|\gamma_{X^{(d)}}(i)\bigr|^{k_i}
\nonumber\\[-8pt]\\[-8pt]
&\leq& \frac{4^h(2h)!}{h!}\Biggl(\sum_{j=0}^d
\sum_{k=0}^{d-j}|\theta_k
\theta_{k+j}|\Biggr)^h \leq\frac{4^h(2h)!}{h!}\Biggl(\sum
_{k=0}^d |\theta_k|
\Biggr)^{2h}.\nonumber
\end{eqnarray}
This
%trivially
implies
%that
(C3) holds,
% $\sum_{h\geq0}\beta_{2h,d}^{-\frac{1}{2h}}=\infty$,
proving the lemma.
%i.e. Carleman's condition is satisfied.
Proof of Theorem \ref{thmautocov}(a) is now complete.\vadjust{\goodbreak}
\end{pf}
%
%s3.1.4 #&#
\subsubsection{Proof of Theorem \texorpdfstring{\protect\ref{thmautocov}}{2.1}\textup{(b)}
(infinite order case)}\label{subsubsectionProofofTheoremautocovb}
%$m_n=n$

\textit{First}, \textit{we assume $\{\varepsilon_t\}$ is i.i.d.} Fix
$\varepsilon>0$. Choose $d$ such that $\sum_{k\geq d+1}|\theta
_k|\leq\varepsilon$.
For convenience we will write
$\Gamma_n(X)=\Gamma_{n}$.
Clearly, $\Gamma_{n}=A_{n}A_{n}^{T}$ where
%Define a matrix $A_n$ such that,
%
\[
(A_{n})_{i,j} = %
\cases{
X_{j-i}, &\quad if $1\leq j-i\leq n$,
\cr
0, &\quad otherwise.}
\]

% $$A_{n}=\frac{1}{\sqrt{n}}\left[ \begin{array} {ccccccccc}
% 0 & X_{1} & X_{2} & \ldots& X_{n-1} & X_{n} & 0 & \ldots& 0 \\
%0 & 0 & X_{1} & \ldots& X_{n-2} & X_{n-1} & X_{n} & \ldots& 0 \\
%& & & & \vdots& & & & \\
%0 & 0 & 0 &\ldots& 0 & X_{1} & X_{2} & \ldots& X_{n}
%so that \begin{align*}(A_{n})_{i,j}&= X_{j-i},\mbox{ if } 1\leq j-i
% &=0,\mbox{ otherwise.}
% and we have
% $$d_{\mathrm{BL}}^{2}(F^{\Gamma_{n,d}},F^{\Gamma_{n}}) \leq\frac{2}{n^2}\Tr
% \Tr\left[(A_{n,d}-A_{n})(A_{n,d}-A_{n})^{T}\right].$$
By \textit{ergodic theorem}, a.s., we have the following two relations:
\begin{eqnarray*}
\frac{1}{n} \bigl[\Tr(\Gamma_{n,d}+ \Gamma_{n})
\bigr] &=& \frac{1}{n} \Biggl[\sum_{t=1}^{n}X_{t,d}^{2}+
\sum_{t=1}^{n}X_t^{2}
\Biggr]\rightarrow\E\bigl[X_{t,d}^{2}+X_{t}^{2}
\bigr] \leq2\sum_{k=0}^{\infty}
\theta_k^2. %&\leq{\frac{2}{n}\sum_{t=1}^{n}(\sum_{k=0}^{\infty}|
%&\leq{\frac{2}{n}\sum_{t=1}^{n}(\sum_{k=0}^{\infty}|
%&= {2(\sum_{k=0}^{\infty}|\theta_k|)\sum_{k=0}^{\infty}|
\\
%
%Hence
%$$
% {\limsup_n \frac{1}{n}(\Tr(\Gamma_{n,d}+ \Gamma_{n}))}
%& \leq{ 2(\sum_{k=0}^{\infty}|\theta_k|)\sum_{k=0}^{
%& \leq{2(\sum_{k=0}^{\infty}|\theta_k|)^2} \mbox{ a.s. }.
%
% Similarly, a.s.,
%
\frac{1}{n}\Tr\bigl[(A_{n,d}-A_{n})
(A_{n,d}-A_{n})^{T}\bigr]&=&\frac{1}{n}\sum
_{t=1}^{n}(X_{t,d}-X_{t})^{2}
\rightarrow\E[X_{t,d}-X_{t}]^2 \leq\sum
_{k=d+1}^{\infty}\theta_k^2
\leq\varepsilon^2.
%&\leq\frac{1}{n} {\sum_{t=1}^{n}(\sum_{k=d+1}^{\infty}|
%&\leq\frac{1}{n} {\sum_{t=1}^{n}(\sum_{k=d+1}^{\infty}|
%&= {(\sum_{k=d+1}^{\infty}|\theta_k|)\sum_{k=d+1}^{
\end{eqnarray*}
%
%Hence
%$$
% {\limsup_n} \frac{1}{n}
%& \leq{(\sum_{k=d+1}^{\infty}|\theta_k|)\sum_{k=d+1}^{
%& \leq{(\sum_{k=d+1}^{\infty}|\theta_k|)^2} \mbox{ a.s. }
%& \leq\varepsilon^2.
Hence using Lemma \ref{lemaux}(b), a.s.
%
%e3.26 #&#
\begin{equation}
\label{eqnautocovgeneral} \limsup_n d_{\mathrm{BL}}^{2}
\bigl(F^{\Gamma_{n,d}},F^{\Gamma
_{n}}\bigr) \leq2\Biggl(\sum
_{k=0}^{\infty}|\theta_k|
\Biggr)^2\varepsilon^2.
\end{equation}
Now $F^{\Gamma_{n,d}}\stackrel{w}{\to}F_d$ a.s.
Since $d_{\mathrm{BL}}$ metrizes weak convergence of
probability measures
%(on complete separable metric spaces, in particular on $\mathbb{R}$)
%we have
as $n\to\infty$, $d_{\mathrm{BL}}(F^{\Gamma
_{n,d}},F_d)\to0$, a.s.
Since $\{F^{\Gamma_{n,d}}\}_{n\geq1}$ is Cauchy
with respect to $d_{\mathrm{BL}}$ a.s.,
by triangle inequality, %$$d_{\mathrm{BL}}(F^{\Gamma_{n}},F^{\Gamma_{m}})\leq
%d_{\mathrm{BL}}(F^{\Gamma_{n}},F^{\Gamma_{n,d}})+d_{\mathrm{BL}}(F^{\Gamma_{n,d}},F^{
and (\ref{eqnautocovgeneral}),
$\limsup_{m,n}d_{\mathrm{BL}}(F^{\Gamma_{n}},F^{\Gamma_{m}})\leq2\sqrt{2}(\sum
_{k=0}^{\infty}|\theta_k|)\varepsilon$.
Hence $\{F^{\Gamma_{n}}\}_{n\geq1}$ is Cauchy with respect to
$d_{\mathrm{BL}}$ a.s. Since $d_{\mathrm{BL}}$ is complete, there exists
a probability measure $F$ on $\mathbb{R}$ such that
$F^{\Gamma_{n}}\stackrel{w}{\to}F$ a.s.
Further
\[
d_{\mathrm{BL}}(F_d,F) = \lim_{n}d_{\mathrm{BL}}
\bigl(F^{\Gamma
_{n,d}},F^{\Gamma_{n}}\bigr)\leq\sqrt{2}\Biggl(\sum
_{k=0}^{\infty}|\theta_k|\Biggr)\varepsilon
\]
and hence $F_d\stackrel{w}{\to}F$ as $d\to
\infty$.
Since $\{F_d\}$ are non-random,
%we conclude that
$F$ is also non-random.

Now if $\{\varepsilon_t\}$ is not i.i.d. but independent and uniformly
bounded by some $C>0$, then the above proof is even simpler.
%One has to simply note that $ {\limsup_n n^{-1}
We omit the details.
%This completes the proof of the first part.

To show convergence of
%the moments
$\{\beta_{h,d}\}$, we note that under
%the
Assumption \ref{assumptionB}(b),
%of summability of the coefficients,
(\ref{eqsupbeta}) yields
%
%e3.27 #&#
\begin{equation}
\label{eqbetasup1} \sup_d |
\beta_{h,d}|\leq c_h:= \frac{4^h(2h)!}{h!}\Biggl(\sum
_{k=0}^\infty|\theta_k|
\Biggr)^{2h}<\infty\qquad \forall h\geq0.
\end{equation}
Hence for every fixed $h$, $\{A_d^h\}$ is uniformly integrable where
%Thus we have uniform integrability of
%all powers of
$A_d\sim F_d$.
% are uniformly integrable.
Since
$F_d\stackrel{w}{\to}F$,
%we thus conclude
%
\[
\beta_h=\int{x^h \mrmd F}=\lim_d
\int{x^h \mrmd F_d}=\lim_{d\to\infty}
\beta_{h,d},
\]
completing the proof of (b). Since $|\beta_h|\leq c_h$, it easily
follows that $\{\beta_h\}_{h\geq0}$ satisfies (C3) and hence uniquely
determines the distribution $F$.

%s3.1.5 #&#
\subsubsection{Proof of Theorem \texorpdfstring{\protect\ref{thmautocov}}{2.1}\textup{(c)}}\label
{subsubsectionProofofTheoremautocovc}
% \subsubsection{Proof of moment ordering}
We first claim that
%prove the following.
for $d\geq0$
%and $k_0,\ldots, k_d\geq0$, we have
$p^{(d)}_{k_0,\ldots, k_d}=p^{(d+1)}_{k_0,\ldots, k_d,0}$.
To see this, consider a graph $G$ with $2h$ vertices with $h$ connected
components and two vertices in each component. Let
\begin{eqnarray*}
\mathcal{M}&=&\bigl\{\mathbf{a}\dvt %=(a_1,...,a_{2h})\in\{1,2,...,2n \}^{2h}
\mathbf{a}\mbox{ is minimal }d\mbox{
matched, induces }G \mbox{ and } |a_x-a_y|=d+1
\\
&&\hspace*{3.5pt}\mbox{for some } x, y \mbox{ belonging to distinct components of }
G\bigr\}.
\end{eqnarray*}
%
%???Drop the parts below and refer to tech report???
Then one can easily argue that $\#\mathcal{M}=\RMO(n^{h-1})$ and
consequently $\#\{(\mathbf{t},\bolds{\pi})\in\mathcal{A}|
\mathbf{a}(\mathbf{t},\bolds{\pi})\in\mathcal{M}\}=\RMO(n^h)$.
Hence,
\begin{eqnarray*}
&&p^{(d)}_{k_0,\ldots, k_d}
\\
&&\quad=\lim_{n \rightarrow\infty}\frac{1}{n^{h+1}}\# \Biggl\{(\mathbf
{t},\bolds{
\pi})\in\mathcal{A}| \mathbf{a}(\mathbf{t},\bolds{\pi})\mbox{ is
minimal }d
\mbox{ matched}\\
&&\quad\hspace*{67.1pt} \mbox{with partition }\{1,\ldots,2h\}=\bigcup
_{l=1}^{h} \{i_l,j_l \}
\\
&&\quad\hspace*{67.1pt}\mbox{and there are exactly } k_s \mbox{ many } l\mbox{'s for which}
\\
&&\quad\hspace*{67.1pt}\bigl|\mathbf{a}(\mathbf{t},\bolds{\pi}) (i_l)-\mathbf{a}(\mathbf
{t},\bolds{\pi}) (j_l)\bigr|=s, s=0,\ldots,d, \I_{\mathbf{a}(\mathbf
{t},\bolds{\pi})}=1\mbox{ and }
\\
&&\quad\hspace*{67.1pt}\bigl|\mathbf{a}(\mathbf{t},\bolds{\pi}) (x)-\mathbf{a}(\mathbf{t},\bolds
{\pi})
(y)\bigr|\geq d+2\mbox{ if }x,y \mbox{ belong to }
\\
&&\quad\hspace*{67.1pt} \mbox{different partition blocks} \Biggr\}
\\
&&\quad=p^{(d+1)}_{k_0,\ldots,k_d,0}.
\end{eqnarray*}

Thus for $\theta_0,\ldots, \theta_d\geq0$ and $d\geq1$,
\begin{eqnarray*}
\beta_{h,d} &\geq&\sum_{S_{h, d-1}}p^{(d)}_{k_0,\ldots,k_{d-1},0}
\prod_{i=0}^{d-1} \bigl[\gamma_{X^{(d)}}(i)
\bigr]^{k_i} \\
&\geq&\sum_{S_{h, d-1}}p^{(d-1)}_{k_0,\ldots,k_{d-1}}
\prod_{i=0}^{d-1}\bigl[\gamma_{X^{(d-1)}}(i)
\bigr]^{k_i} =\beta_{h,d-1},
\end{eqnarray*}
proving the result.

%proving the result.
%moments)
Incidentally, if Assumption \ref{assumptionB}(a) is violated, then the
%moment
ordering need not hold. This can be checked by considering an MA(2) and
an MA(1) process with parameters $\theta_0, \theta_1, \theta_2$ and
%with parameter set $\theta_0, \theta_1$
where
$\theta_2 = - \kappa\theta_0$, $\theta_0, \theta_1 >0$. Then
%it can be shown that
$ \beta_{2,2} < \beta_{2,1}$ if we choose $\kappa> 0$ sufficiently
small. The details are
%which is
available in Basak, Bose and Sen \cite{Basakbosesen}.
%
% Using Lemma \ref{lemcoeffmad} we note that
%p_{0,2,0}^{(2)}( \theta_0 \theta_1 + \theta_1 \theta_2)^2 +
%p_{0,0,2}^{(2)} \theta_0^2 \theta_2^2
%and
%Using Lemma \ref{lemiteration} we get
%$p_{0,2,0}^{(2)}=p_{0,2}^{(1)}$. Further, it is also not hard to
%verify that $p_{0,0,2}^{(2)}=p_{0,2,0}^{(2)}$. Thus $\beta_{2,2} \ge
%p_{2,0}^{(1)} [\theta_2^4 + 2 (\theta_0^2+\theta_1^2) \theta_0^2] +
%p_{0,2}^{(1)} [\theta_1^2 \theta_2^2 + 2 \theta_0 \theta_1^2 \theta_2
%+ \theta_0^2 \theta_2^2] \ge0
%Now taking $\theta_2 = - \kappa\theta_0$ where $\kappa>0$ and $
%p_{0,2}^{(1)} (\kappa-2)) \theta_1^2 \ge0
%After solving a linear equation for $\kappa$, it is easily seen that
%there exists a $\kappa^*>0$ such that if $\kappa\in(0, \kappa^*)$,
%then coefficient of $\theta_1$ will be negative. Hence fixing some
%arbitrary value of $\theta_0>0$, and $\kappa\in(0,\kappa^*)$ one can
%increase the value of $\theta_1$ arbitrarily to get $ \beta_{2,2} <
%$d=1,2$ has made calculation much simpler. I still do not know whether
%it is true (or false) that for high enough $h,d$ the ordering of
%moments will be valid.

%s3.1.6 #&#
\subsubsection{Proof of unbounded support of $F_d$ and $F$}
\label{subsubsectionProofofunboundedsupport}

For any word $w$, let $|w|$ denote the length of the word. Let
\begin{eqnarray*}
\mathcal{W}&=&\bigl\{w=w_1w_2\dvt |w_1|=2h=|w_2|;\\
&&\hspace*{4.85pt}
w,w_1,w_2\mbox{ are zero pair matched};
w_1[x] \mbox{ matches }
\\
&&\hspace*{4.85pt} \mbox{with } w_1[y]\mbox{ iff } w_2[x] \mbox{ matches with
}w_2[y]\bigr\}.
\end{eqnarray*}
Then
%
%e3.28 #&#
\begin{equation}
\label{eqnunbounded} \beta_{2h,d}\geq\bigl[\gamma_{X^{(d)}}(0)
\bigr]^{2h}p_{2h,0,\ldots, 0}\geq\bigl[\gamma_{X^{(d)}}(0)
\bigr]^{2h}\sum_{w\in\mathcal{W}}\lim
_n n^{-(2h+1)}\#\Pi^\ast(w).
\end{equation}

For $w=w_1 w_2\in\mathcal{W}$, let $\{1,\ldots, 2h\}
=\bigcup_{i=1}^{h}(i_s,j_s)$ be the partition corresponding to $w_1$. Then
%
%e3.29 #&#
\begin{eqnarray*}
&&\lim_n\frac{\#\Pi^\ast(w)}{n^{2h+1}}\geq\lim
_n\frac
{1}{n^{2h+1}}\#\bigl\{(\mathbf{t}, \bolds{\pi})\dvt
t_{i_s}=t_{j_s}\mbox{ and } \pi_{i_s}-
\pi_{i_s-1}=\pi_{j_s-1}-\pi_{j_s}\\
&&\hspace*{111pt}\mbox{ for } 1\leq s
\leq h;
t_j+|\pi_{j}-\pi_{j-1}|\leq n\mbox{, for }
1\leq j\leq2h\bigr\}.
\end{eqnarray*}
Now adapting the ideas of Bryc, Dembo and Jiang \cite{bry}, we obtain
that for each $d$ finite $F_d$ has unbounded support. Since $\{\beta
_{h,d}\}$ increases to $\beta_h$, same conclusion is true for $F$. For
details see Basak, Bose and Sen \cite{Basakbosesen}.
\subsection{Outline of the proof of Theorem \texorpdfstring{\protect\ref{thmautocovband}}{2.3}}
\label{subsectionOutlineoftheproofwhenalpha<1}

%Here is an outline of the changes required in the proof of
%Theorem \ref{thmautocov}.
% for other values of $\alpha$.

%s3.2.1 #&#
\subsubsection{Proof of Theorem \texorpdfstring{\protect\ref{thmautocovband}}{2.3}\textup{(a)},
\textup{(b)} for the case \texorpdfstring{$0<\alpha<1$}{0<alpha<1}}

Let $\beta_h(\Gamma_{n,d}^{\alpha,I})$ and $\beta_h(\Gamma
_{n,d}^{\alpha,\mathit{II}})$ be the $h$th moments, respectively, of the ESD
of type I and type II ACVMs
%autocovariance matrices
with parameter $\alpha$. We begin by noting that the expression for
these contain an extra indicator term $\I_1=\prod_{i=1}^{h} \I(|\pi
_{i-1} - \pi_i| \le m_n)$ and $\I
_2=\prod_{i=1}^{h} \I( 1 \le\pi_i \le m_n)$,
respectively. For type II ACVMs
%autocovariance matrices
since there are $m_n$ eigenvalues instead of $n$, the normalising
denominator is now~$m_n$.
Hence,
\[
\beta_h\bigl(\Gamma_{n,d}^{\alpha,I}
\bigr) = \frac{1}{n^{h+1}} \mathop{\sum_{1\leq\pi_0,\ldots,\pi_h
\leq n }}_{ \pi_h=\pi_0}
\Biggl[\prod_{j=1}^h \Biggl(\sum
_{t_j=1}^n X_{t_j,d}X_{t_j+|\pi_j-\pi
_{j-1}|,d}
\I_{(t_j+|\pi_j-\pi_{j-1}|\leq n)} \Biggr) \Biggr]\I_1
\]
and
\[
\frac{m_n}{n}\beta_h\bigl(\Gamma_{n,d}^{\alpha,\mathit{II}}
\bigr) = \frac{1}{n^{h+1}}\mathop{\sum_{1\leq\pi_0,\ldots,\pi_h \leq n
}}_{\pi_h=\pi_0}
\Biggl[\prod_{j=1}^h \Biggl(\sum
_{t_j=1}^n X_{t_j,d}X_{t_j+|\pi_j-\pi
_{j-1}|,d}
\I_{(t_j+|\pi_j-\pi_{j-1}|\leq n)} \Biggr) \Biggr]\I_2.
%& & \times\prod_{i=1}^{h} \I( 1 \le\pi_i \le m_n)
\]
%
%Since $\alpha\ne0$ for $\Gamma_{n}^{\alpha, \mathit{II}}(X^{(d)})$
It is thus enough to establish the limits on the right side of the
above expressions.
%of $\alpha_n \beta_h(\Gamma_{n,d}^{\alpha,\mathit{II}})$ instead of $ \beta_h(
and we can follow similar steps as in the proof of Theorem \ref{thmautocov}.

Since there are only the
%some
extra indicator terms, the negligibility of higher order edges and
verification of (C2) and (C3) needs no new arguments. Likewise,
verification of (C1) is also similar except that there is now an extra
indicator term in the expression for $p_w^{(d)}$. %From (\ref{eqtype1})
%and (\ref{eqtype2})
%we note that the extra term ${\prod_{j=2}^{h+1} \I(|
%will be inside the expectation for $\Gamma_{n,d}^{\alpha, I}$ and for $
%{\prod_{j=1}^{h+1} \I(\bolds{\lambda^b_{j+h}} \in(0,
This takes care of the finite $d$ case.
For $d=\infty$, note that the type II ACVMs
%autocovariance matrices
are $m_n \times m_n$ principal subminor of the original sample ACVMs
%autocovariance matrices
and hence are automatically non-negative definite. We can write $\Gamma
_{n}^{\alpha, \mathit{II}}(X^{(d)})= (A_{n,d}^{\alpha,\mathit{II}})(A_{n,d}^{\alpha,\mathit{II}})^T$
where $A_{n,d}^{\alpha, \mathit{II}}$ is the first $m_n$ rows of $A_{n,d}$. Thus
imitating the proof of Theorem \ref{thmautocov}, we can move from
%limiting distribution
finite $d$ to $d= \infty$. However for type I
ACVMs,
%autocovariance matrices
we cannot apply these arguments,
%heorem \ref{thmautocov}
as these matrices are
%may
not necessarily
non-negative definite. Rather we proceed as in the proof of Theorem \ref
{thmautocov2}.
Previous proof of unbounded support now needs only minor changes. We
omit the details.

Since $\Gamma_{n,d}^{\alpha, \mathit{II}}$ is non-negative definite, the
technique of proof of Theorem \ref{thmautocov} can
%may
be adopted under Assumption \ref{assumptionA}(a).

%(i) ????? Here one can again write
%$$\Gamma_n^{\alpha}(X^{(d)}) =(A_{n,d}^\alpha)( A_{n,d}^{\alpha})^T$$
%where $A_{n,d}^\alpha$ is same as $A_{n,d}$ except the last $n-m_n$
%rows which will now be identically $0$.?????

%(ii) When we write the expression for the $h^{th}$ moment, now there
%the only difference than the usual sample autocovariance case is
%presence of
%an extra indicator term ${\prod_{i=1}^{h} \I(|\pi_{i-1} -

%(iii) As a consequence, the negligibility of higher order edges,
%verification of (C2) and (C3) needs no new arguments.

%(iv) Verification of (C1) is also similar. There will be an extra term
%coming in the expression for $p_w^{(d)}$ now. The extra term $
%{\prod_{j=2}^{h+1} \I(|\bolds{\lambda^b_{j+h-1}} -

%(v) Proof of the other parts are either same as before or need only
%minor modifications.

%s3.2.2 #&#
\subsubsection{Proof of Theorem \texorpdfstring{\protect\ref{thmautocovband}}{2.3}\textup{(b)} for
%for $\alpha=0$,
type \textup{I} band ACVM}
%autocovariance matrix}
\label{subsectionaplha0proof}

\textit{Existence}: Let
%Denoting
$p_w^{(d),0,I}$ be the limiting contribution of the word $w$ for type I ACVM
%autocovariance matrix
with band parameter $\alpha=0$. Then
\[
p_w^{(d),0,I}:=\lim_n
\frac{1}{n^{h+1}} \sum_{\mathbf{b}\in{\cal B}(w)}\E_{U_n} \bigl[
\I_{m,h}(U_n) \times\I_{\bolds{\lambda}^{\mathbf{b}}, m}(U_n)
\times\I_{\bolds{\lambda}^{\mathbf{b}},L}(U_n) \times\I_{{\cal
T}}^{I}(U_n)
\bigr],
\]
where
\[
\I_{{\cal T}}^I(U_n)=\I_{{\cal T},L}(U_n)
\times\I_{{\cal T},m}:=\mathop{\prod_{j=1}}_{j \notin{\cal T}}^{h}
\I\bigl(0 \le L_j(U_n) \le m_n\bigr) \times
\prod_{j \in{\cal T}} \I( -m_n \le n_j
\leq0).
\]
If $w$, $\bolds{\lambda}^{\mathbf{b}}_{j+h-1} \ne\bolds{\lambda}^{\mathbf{b}}_{j+h}$
for some $j$, then $\I_{{\cal T},L}( U_n ) \rightarrow0$ as $n \to
\infty$ and thus limiting contribution from that word will be $0$.
Thus, only those words $w$ for which $\bolds{\lambda}^{\mathbf{b}}_{h+1} =
\bolds{\lambda}^{\mathbf{b}}_{j+h}$ for all $j \in\{1,2,\ldots,h+1\}$ may
contribute non-zero quantity in the limit. This condition also implies
that, for such words no $\pi_i$ belongs to the generating set except
$\pi_0$. This observation together with Lemma 6 of Basak, Bose and
Sen \cite{Basakbosesen}, and the expression for limiting moments
for $\Gamma_n(X)$ shows that $w \in{\cal W}_0^h$ may contribute
non-zero quantity, where
\[
{\cal W}_0^h= \bigl\{w\dvt  |w|=2h, w[i] \mbox{
matches with } w[i+h], n_i \le0, i =1,2,\ldots,h\bigr\}.
\]
Further note that if $w \in{\cal W}_0^h$ then $\mathcal{T}=\{
h+1,h+2,\ldots,2h\}$, and thus $\I_{{\cal T},L} \equiv1$.

For $d=0$ note that $\#{\cal W}_0^h=1$ for every $h$ and one can easily
check that the contribution from that word is $1$. Thus $\beta_{h,0}^0
=\theta_0^{2h}$ and as a consequence, the LSD is $\delta_{\theta_0^2}$.

Now let us consider any $ 0 < d < \infty$. Note that for any $d$
finite, and if $m_n \ge d$, then
\[
\I_{\bolds{\lambda}^{\mathbf{b}},m} \times\I_{\bolds{\lambda
}^{\mathbf{b}},L} \times\I_{{\cal T},m} \rightarrow
\prod_{j=1}^h \I(n_j \le0)\qquad\mbox{as } n \rightarrow\infty.
\]
Combining the above arguments we get that for any $w \in{\cal W}_0^h$,
$p_w^{(d),0,I}$ is the number of choices of $\mathbf{b} \in{\cal
B}(w)$, and $\{n_1,n_2,\ldots,n_h; n_i \le0\}$, such that $\sum_i
n_i b_i=0$.
%Now we argue that the limiting distribution has bounded support. For
%this, it is enough to argue that $\lim_{h \to\infty} (
%that
%k_1,\ldots, k_d}} \Bigg)
%& = & 2^{2h} \Big( \sum_{i=0}^d |\gamma_{X^{(d)}}(i)| \Big)^{2h} \le
%2^{2h} \Big(\sum_{i=0}^d\sum_{j=0}^{d-i} |\theta_i \theta_{i+j}|
%
%Thus
%$$\limsup_{h \to\infty} (\beta_{2h,d}^0)^{1/2h} \le2(\sum_{i=0}^d |
%All the above conclusions can be proved even if $m_n \to d$ by a
%little modification of the above ideas.
%Since $b_i \in\{-1,1\}$ and $|n_i| \in\{0,1,\ldots,d\}$ clearly $p_{

Noting that type I ACVMs
%autocovariance matrices
are not necessarily non-negative definite, we need to adapt the proof of
Theorem \ref{thmautocov2}. Details are omitted.
%, existence of LSD for $d=\infty$, and it's connection with finite $d$
%is established.

%For $d=\infty$ also the support is bounded. To argue this, one simply
%notes that (\ref{eq0dbound}) will also imply
%and hence we have the claim.

\textit{Identification of the LSD}: Now it remains to argue that the
limit we obtained is same as $f_X(U)$.
%the limiting spectral distribution of $\Sigma_n$.
For $d=0$ LSD is $\delta_{\theta_0^2}$
%autocovariance matrices
and it is trivial to check it is same as $f_X(U)$.

For $0 < d < \infty$, note that the proof does not use the fact that
$m_n \to\infty$ and we further note that for any sequence $\{m_n\}$
the limit we obtained above will be same whenever $\liminf_{n \to
\infty} m_n \ge d$. So in particular the limit will be same if we
choose another sequence $\{m_n'\}$ such that $m_n'=d$ for all $n$. Let
$\Gamma_{n',d}^{I}$ denote the type I
ACVM
%autocovariance matrix
where we put $0$ instead of $\hat{\gamma}_{X^{(d)}}(k)$ whenever $k
>m_n'$ and let $\Sigma_{n,d}$ be the $n \times n$ matrix whose
$(i,j)$th entry is the population autocovariance $\gamma
_{X^{(d)}}(|i-j|)$. Now from Lemma \ref{lemaux}(a), we get
\begin{eqnarray*}
d^2_{\mathrm{BL}}\bigl(F^{\Gamma_{n',d}^{I}}, F^{\Sigma_{n,d}}\bigr)
&\le& \frac{1}{n} \Tr\bigl(\Gamma_{n',d}^{I} -
\Sigma_{n,d}\bigr)^2
\\
& \le& 2 \bigl(\hat{\gamma}_{X^{(d)}}(0)- \gamma_{X^{(d)}}(0)
\bigr)^2+ \cdots+ 2 \bigl(\hat{\gamma}_{X^{(d)}}(d)-
\gamma_{X^{(d)}}(d)\bigr)^2.
\end{eqnarray*}
For any $j$ as $n \to\infty$, $\hat{\gamma}_{X^{(d)}}(j) \to\gamma
_{X^{(d)}}(j)$ a.s.
%Using the fact that
Since $d$ is finite, the right
%hand
side of the above expression goes to $0$ a.s. This proves the claim
for $d$ finite.

To
prove
the result for the case $d=\infty$, first note that we already have
\[
\operatorname{LSD}\bigl(\Gamma_{n,d}^{0,I}\bigr)=\operatorname{LSD}(\Sigma_{n,d}):=
G_d \quad\mbox{and}\quad \operatorname{LSD}\bigl(\Gamma_{n,d}^{0,I}\bigr)
\stackrel{w} {\to} \operatorname{LSD}\bigl(\Gamma_{n}^{0,I}\bigr)
\qquad\mbox{as } d \to\infty.
\]
Thus, it is enough to prove that
$G_d \stackrel{w}{\to} G(=\operatorname{LSD}(\Sigma_n))$ as $d \to
\infty$
where $\Sigma_n$ is the $n \times n$ matrix whose $(i,j)$th entry
is $\gamma_{X}(|i-j|)$.
Define a sequence of $n \times n$ matrices $\bar{\Sigma}_{n,d}$ whose
$(i,j)$th entry is $\gamma_{X}(|i-j|)$ if $|i-j| \le d$ and
otherwise $0$. By triangle inequality,
\[
d^2_{\mathrm{BL}} \bigl(F^{\Sigma_{n,d}},
F^{\Sigma_n}\bigr) \le2 d^2_{\mathrm{BL}}\bigl(F^{\Sigma
_{n,d}},
F^{\bar{\Sigma}_{n,d}}\bigr) + 2 d^2_{\mathrm{BL}} \bigl(F^{\bar{\Sigma
}_{n,d}},
F^{\Sigma_n}\bigr).
\]
Fix any $\varepsilon>0$. Fix $d_0$ such that $ (
\sum_{j=0}^\infty|\theta_j| )^2 ( \sum_{l=d+1}^\infty
|\theta_l| )^2 \le\frac{\varepsilon^2}{32}$ for all $d \ge d_0$.
Now again using Lemma \ref{lemaux}(a) we get the following two relations:
\begin{eqnarray*}
\limsup_{n} d^2_{\mathrm{BL}}
\bigl(F^{\Sigma_{n,d}}, F^{\bar{\Sigma}_{n,d}}\bigr) %& \le& \limsup_n
& \le& 2
\bigl[ \bigl(\gamma_{X^{(d)}}(0) - \gamma_X(0)
\bigr)^2 + \cdots+ \bigl(\gamma_{X^{(d)}}(d) -
\gamma_X(d)\bigr)^2 \bigr]
\\
& = & 2 \sum_{j=0}^d
\Biggl(\sum_{k=d-j+1}^\infty\theta_k
\theta_{j+k} \Biggr)^2 \le\frac{\varepsilon^2}{16},
\\
%& \le& 2 \sum_{j=0}^d \Big(\sum_{k=d-j+1}^\infty|\theta_k| |
%& \le& 2 \Big(\sum_{j=0}^d \sum_{k=d-j+1}^\infty|\theta_k| |
%& \le& 2 \Big( \sum_{j=0}^d \sum_{l=d+1}^\infty|\theta_{l-j}| |
%& = & 2 \Big( \sum_{l=d+1}^\infty|\theta_{l}| (|\theta_{l}|+ \cdots+ |
%& \le& 2 \Big( \sum_{l=d+1}^\infty|\theta_{l}| \sum_{j=0}^\infty|
%
%and similarly
d^2_{\mathrm{BL}}
\bigl(F^{\bar{\Sigma}_{n,d}}, F^{\Sigma_n}\bigr) & \le& \limsup
_n \frac{1}{n} \Tr(\bar{\Sigma}_{n,d} -
\Sigma_n)^2 \le\frac
{\varepsilon^2}{16}.
%& \le& 2 \sum_{j=d+1}^ \infty\gamma_X(j)^2
%& = &2 \sum_{j=d+1}^\infty\Big[ \sum_{k=0}^\infty\theta_k
%& \le& 2 \sum_{j=d+1}^\infty\Big[ \sum_{k=0}^\infty|\theta_k| |
%& \le& 2 \Big[ \sum_{j=d+1}^\infty\sum_{k=0}^\infty|\theta_k| |
%& = & 2 \Big[\sum_{k=0}^\infty|\theta_k| \sum_{j=d+1}^\infty|
%& \le& 2 \Big[\sum_{k=0}^\infty|\theta_k| \sum_{j=d+1}^\infty|
\end{eqnarray*}
Thus, %\begin{equation} \label{eqorg4}
$\limsup_n d_{\mathrm{BL}} (F^{\Sigma_{n,d}}, F^{\Sigma_n})
\le\varepsilon/2$, for any $d \ge d_0$, and therefore by triangle
inequality, $d_{\mathrm{BL}}(F^{G_d},F^G) \le\varepsilon$.
This completes the proof.

%s3.2.3 #&#
\subsubsection{Proof of Theorem
\texorpdfstring{\protect\ref{thmautocovband}}{2.3}\textup{(b)} for type \textup{II} band autocovariance
matrix}\label{subsectionalphazerotypeii} %Part of this proof will be
%different from the corresponding proof for type I band ACVMs,
%autocovariance matrices,
%as the expressions for the $h^{th}$ moment
%of the ESD of these matrices differ by a factor $\alpha_n$ and here
%$\alpha_n \rightarrow0$.

First, note that by Lemma \ref{lemorder-reduction} we need to
consider only minimal matched terms. Let
%Let $G_t$ be the set of $t_i$ belonging to the generating vertices and
%let $G_\pi$ be the same for the $\pi_i$'s. Thus formally
%
\[
G_t=\{t_i\dvt  t_i \in G\} \quad\mbox{and}\quad
G_{\pi} = \{ \pi_i\dvt  \pi_i \in G\}.
\]
Since $1 \le\pi_i \le m_n$ for all $i$, by similar arguments as in
Lemma \ref{lemorder-reduction} we get
\[
\mbox{number of choices of } \mathbf{a}(\mathbf{t}, \bolds\pi)= \RMO\bigl(n^{\#G_t}
m_n^{\#G_\pi}\bigr).
\]
Thus, for any word $w$ such that $\#G_t < h$ the limiting contribution
will be $0$. Hence only contributing words e in this case are those for
which $\#S_3(w)=\#S_4(w)=h$. and from Lemma 6 of Basak, Bose and
Sen \cite{Basakbosesen}, the only contributing words are those
belonging to $\mathcal{W}_0^h$. Therefore using same arguments as in
the proof of Theorem \ref{thmautocovband}, for type I ACVM, for
$\alpha=0$ we obtain the same limit. All the remaining conclusions
here follow from the proof for type I ACVMs
%autocovariance matrices
with parameter $\alpha=0$. %But here note that in order to achieve the
%same limit as obtained for type I ACVMs
%autocovariance matrices
%we need $m_n \to\infty$ even if $d$ is finite.

Since type II ACVMs are non-negative definite, connection between the LSD
%limiting distribution
for finite $d$ and $d= \infty$ is proved adapting the ideas from the
proof of Theorem \ref{thmautocov}.

%s3.2.4 #&#
\subsubsection{Proof of Theorem \texorpdfstring{\protect\ref{thmautocovband}}{2.3}\textup{(c)}}
Since $K$ is bounded, negligibility of higher order edges and
verification of (C2) and (C3) is same as before. Verification of
(C1) is also same, with an extra indicator in the limiting
expression. Denoting $p_w^{(d),K}$ to be the limiting contribution from
a word $w$, we have,
\[
p_w^{(d),K} =\lim_n
\E_{U_n} \bigl[ \I_{m,h}(U_n) \times\I
_{\bolds{\lambda}^{\mathbf{b}}, m}(U_n) \times\I_{\bolds
{\lambda}^{\mathbf{b}},L}(U_n)
\times\I_{{\cal T}}^{\I}(U_n) \times
\I_K(U_n) \bigr],
\]
where
\[
\I_K(U_n):= \prod_{j=1}^h
K \biggl(\frac{L_j(U_n)}{m_n} \biggr).
\]
Since $m_n \to\infty$, and $K(\cdot)$ is continuous at $0$,
$K(0)=1$, note that $\I_K \rightarrow1$. Now arguing as in
Section \ref{subsectionaplha0proof}, we get
$p_w^{(d),0,I}=p_w^{(d),K}$ for every word $w$ and thus the limiting
distributions are same in both the cases.
For the case $d=\infty$ the arguments are similar as in Section \ref
{subsectionaplha0proof} and the details are omitted.

\subsection{Proof of Theorem \texorpdfstring{\protect\ref{thmautocov2}}{2.2}}
\label{subsectionProofofTheoremautocov2}

Proceeding as earlier it is easy to see the limit exists, and for each
word $w$, the limiting contribution is given by,
\[
p_w^{*,(d)} = \sum_{\mathbf{b} \in{\cal B}(w)}
\E_U \bigl[ \I_{m,h}(U) \times\I_{\bolds{\lambda}^{\mathbf{b}}} (U)
\times\tilde{I}_{{\cal T}}(U) \bigr].
\]
Comparing the above expression with the corresponding expression for
the sequence $\Gamma_{n,d}$,
\[
\beta_{h,d}\leq\beta_{h,d}^\ast\qquad\mbox{if }
\theta_j\geq0, 0\leq j\leq d.
\]
Relation (\ref{eqsupbeta}) holds with $\beta_{h,d}$ replaced by
$\beta_{h,d}^\ast$. We can use this to prove tightness of $\{
F_{d}^{\ast}\}$ under Assumption \ref{assumptionB}(a)
and thus also Carleman's condition is satisfied.

Since $\Gamma_n^{*}$ and $\Gamma_{n,d}^{*}$ are no longer positive
definite matrices the ideas used in the proof of Theorem \ref{thmautocov}(b) cannot be
adapted here. We proceed as follows instead: Note that
\[
\E\bigl[\beta_{h}\bigl(\Gamma_{n}^{\ast}\bigr)
\bigr]= \frac{1}{n^{h+1}}\E\Biggl[\sum_{(\mathbf{t},\bolds
{\pi})\in\mathcal{A}} \prod
_{j=1}^{h}X_{t_j}\prod
_{j=1}^{h}X_{t_j+|\pi_{j-1}-\pi_{j}|} \Biggr].
\]
Write
\[
X_{t_j}=\sum_{k_j\geq0}\theta_{k_j}
\varepsilon_{t_j-k_j} \quad\mbox{and}\quad X_{t_j+|\pi_{j-1}-\pi_{j}|}=\sum
_{k^{\prime}_j\geq
0}\theta_{k^{\prime}_j}\varepsilon_{t_j+|\pi_{j-1}-\pi_{j}|-k^{\prime}_j}.
\]
Then using the absolute summability Assumption \ref{assumptionB}(b) and applying DCT,
we get
\[
\E\bigl[\beta_{h}\bigl(\Gamma_{n}^{\ast}\bigr)
\bigr]= \mathop{\sum_{k_j,k^{\prime}_j\geq0}}_{ j=1,\ldots,h}\prod
_{j=1}^{h}(\theta_{k_j}
\theta_{k^{\prime}_j})\frac{1}{n^{h+1}} \E\Biggl[\sum
_{(\mathbf{t},\bolds{\pi})\in\mathcal{A}} \prod_{j=1}^{h}
\varepsilon_{t_j-k_j}\varepsilon_{t_j+|\pi_j-\pi
_{j-1}|-k^{\prime}_j} \Biggr]. %
\]
Using the fact that $\{\varepsilon_t\}_{t=1}^\infty$ are uniformly
bounded and absolute summability of $\{\theta_k\}_{k=1}^\infty$ we
note that it is enough to show that the limit below exists.
\[
\lim_n n^{-(h+1)} \E\Biggl[\sum
_{(\mathbf{t},\bolds{\pi})\in\mathcal{A}} \prod_{j=1}^{h}(
\varepsilon_{t_j-k_j}\varepsilon_{t_j+|\pi_j-\pi
_{j-1}|-k^{\prime}_j}) \Biggr].
\]
One can proceed as in the proof of Theorem \ref{thmautocov} to show
that only pair matched words contribute and hence enough to argue that
$\lim n^{-(h+1)}\#\{(\mathbf{t},\bolds{\pi})\in\mathcal{A}\dvt  \{
t_j-k_j, t_j+|\pi_j-\pi_{j-1}|-k^{\prime}_j, j=1,\ldots,h\}$ is
pair matched$\}$ exists, and which follows by adapting the ideas used
in the proof of Theorem \ref{thmautocov}. Note that appropriate
compatibility is needed
among $\{k_j,k_j^{\prime}, j=1,\ldots,h\}$, the word $w$ and the signs
$b_i$ $(=\pm1)$ to ensure that the condition $\pi_0=\pi_h$ is satisfied.
So the above limit will depend on $\{k_j,k_j^{\prime}, j=1,\ldots,h\}$.

We also note that
\begin{eqnarray*}
&&
\lim_n \frac{1}{n^{h+1}} \mathop{\sum
_{w \ \mathrm{pair}\ \mathrm{matched,}}}_{|w|=2h}\#\bigl\{ (\mathbf{t},\bolds{\pi
})\in
\mathcal{A}\dvt \bigl(t_j-k_j, t_j+|\pi
_j-\pi_{j-1}|-k^{\prime}_j
\bigr)_{j=1,\ldots,h}\in\Pi(w)\bigr\}
\\
&&\quad\leq \frac{4^h (2h)!}{h!}.
\end{eqnarray*}
Hence, $F^{*}$ is uniquely determined by its
moments and using DCT, $\beta_{h,d}^{*} \rightarrow
\beta_h^{*}$. Whence it also follows that
$F^{\ast}_d\stackrel{w}{\to}F^{\ast}$. Proof of part (c) is similar
to the proof of Theorem \ref{thmautocov}(c). \qed
%
%re3.1 #&#
\begin{rem}\label{remautocov2}
% We have not proved
Theorem \ref{thmautocov2} has not been proved under Assumption
\ref{assumptionA}\textup{(a)}
% This is
because there is
%now
no straightforward way to apply
(\ref{blmetric})
%. We cannot apply
or (\ref{blmetric2}) since $\Gamma_n^*(X)$ is not non-negative definite.
Simulation results indicate that the same LSD continues to hold under
Assumption \ref{assumptionA}\textup{(a)}.
%for the i.i.d. finite second moment case.
%We conjecture that the LSDs $F_d^\ast$ and $F^\ast$ exist
%under this condition. Moreover, $F_d^\ast$ should converge weakly
%to $F^\ast$
% under this condition.
\end{rem}
\section*{Acknowledgements}

We thank Dimitris Politis and Mohsen Pourahmadi for sharing their work
and thoughts. The constructive comments of the four Referees and the
Associate Editor is gratefully acknowledged. We thank the Editor for
his encouragement to submit a revision.

A. Basak supported by Melvin and Joan Lane endowed Stanford Graduate Fellowship
fund.
A. Bose's research supported by J.C. Bose Fellowship, Govt. of India.
S. Sen supported by NYU graduate fellowship under Henry M. MacCracken Program.

\begin{supplement}%[id=suppA]
\stitle{Simulations}
\slink[doi]{10.3150/13-BEJ520SUPP} %[doi,text={...}] - jei reikia
%suskaldyti doi
\sdatatype{.pdf}
\sfilename{BEJ520\_supp.pdf}
\sdescription{Recall that none of the LSDs have a nice description.
Following the
suggestion of one of the Referees, we have collected some simulation
results in a supplementary file Basak, Bose and Sen \cite
{basakbosesensim}.\\
\indent The simulations are for the AR(1) and MA(1) models. These simulations
provide evidence that the limits are indeed universal and exhibit some
mass on the negative axis for the ESD (and hence the LSD) of
$\Gamma_n^*(X)$. They also show how the LSD of type I banded
$\Gamma_n(X)$ changes with the model as well as the value of the
parameter $\alpha$. The unbounded nature of the LSD is also evident
from these simulations.\\
\indent
For the banded matrices, the simulations demonstrate that for small
values of $\alpha$, the LSD of $\Sigma_n(X)$ and $\Gamma_n(X)$ are
virtually indistinguishable for large $n$, confirming that thinly
banded ACVMs are consistent for $\Sigma_n(X)$. As the value of
$\alpha$ increases, the right tail of the LSD thickens, and the
probability of being near zero decreases. In general, there may be
considerable amount of mass in the negative axis. This mass reduces as
the value of $\alpha$ decreases.\\
\indent
The LSD of $\Gamma_n (X)$ varies as the parameter of the models change.
For both AR(1) and MA(1) models, as $\theta$ increases from $0$, the
tail thickens, and the mass near zero decreases. For the AR(1) model,
when $\theta$ approaches $1$, that is, when the process is near
non-stationary the LSD becomes very flat, and its tail becomes huge.}
\end{supplement}

% imsref loaded by lrinkeviciute, 2013-06-10 14:50:39
%
% imsref loaded by lrinkeviciute, 2013-06-10 20:32:30
% imsref loaded by lrinkeviciute, 2013-06-13 13:37:55
% imsref loaded by lrinkeviciute, 2013-06-13 13:38:08

\printhistory


\begin{thebibliography}{22}
% BibTex style file: bej.bst, 2012-09-27
% Default style options (sort=1,type=number).
% Used options (sort=1,type=number).

%b1 ###
\bibitem{Arcones}
\begin{barticle}[mr]
\bauthor{\bsnm{Arcones},~\bfnm{Miguel~A.}\binits{M.A.}}
(\byear{2000}).
\btitle{Distributional limit theorems over a stationary {G}aussian sequence of
  random vectors}.
\bjournal{Stochastic Process. Appl.}
\bvolume{88}
\bpages{135--159}.
\bid{doi={10.1016/S0304-4149(99)00122-2}, issn={0304-4149}, mr={1761993}}
\bptok{imsref}%
\end{barticle}
\endbibitem

%b2 ###
\bibitem{Baizhou2008}
\begin{barticle}[mr]
\bauthor{\bsnm{Bai},~\bfnm{Zhidong}\binits{Z.}} \AND
  \bauthor{\bsnm{Zhou},~\bfnm{Wang}\binits{W.}}
(\byear{2008}).
\btitle{Large sample covariance matrices without independence structures in
  columns}.
\bjournal{Statist. Sinica}
\bvolume{18}
\bpages{425--442}.
\bid{issn={1017-0405}, mr={2411613}}
\bptok{imsref}%
\end{barticle}
\endbibitem

%b3 ###
\bibitem{Bai99}
\begin{barticle}[mr]
\bauthor{\bsnm{Bai},~\bfnm{Z.~D.}\binits{Z.D.}}
(\byear{1999}).
\btitle{Methodologies in spectral analysis of large-dimensional random
  matrices, a review (with discussions)}.
\bjournal{Statist. Sinica}
\bvolume{9}
\bpages{611--677}.
\bid{issn={1017-0405}, mr={1711663}}
\bptok{imsref}%
\end{barticle}
\endbibitem

%b4 ###
\bibitem{Basak}
\begin{bmisc}[auto:STB|2013/06/05|13:45:01]
\bauthor{\bsnm{Basak},~\bfnm{Anirban}\binits{A.}}
(\byear{2009}).
\btitle{Large dimensional random matrices. M. Stat. Project report, May 2009.
  Indian Statistical Institute}.
\bptok{imsref}%
\end{bmisc}
\endbibitem

%b5 ###
%  \bauthor{\bsnm{Bose},~\bfnm{Arup}\binits{A.}}
%(\byear{2010}).

%b6 ###
\bibitem{Basakbosesen}
\begin{bmisc}[auto:STB|2013/06/05|13:45:01]
\bauthor{\bsnm{Basak},~\bfnm{Anirban}\binits{A.}},
  \bauthor{\bsnm{Bose},~\bfnm{Arup}\binits{A.}} \AND
  \bauthor{\bsnm{Sen},~\bfnm{S.}\binits{S.}}
(\byear{2011}).
\bhowpublished{Limiting spectral distribution of sample autocovariance
  matrices. Technical Report R11 2011. Stat-Math Unit, Indian Statistical
  Institute. Available at \url{http://arxiv.org/pdf/1108.3147v1.pdf}}.
\bptok{imsref}%
\end{bmisc}
\endbibitem

%b7 ###
\bibitem{basakbosesensim}
\begin{bmisc}[auto:STB|2013/06/05|13:45:01]
\bauthor{\bsnm{Basak},~\bfnm{Anirban}\binits{A.}},
  \bauthor{\bsnm{Bose},~\bfnm{Arup}\binits{A.}} \AND
  \bauthor{\bsnm{Sen},~\bfnm{S.}\binits{S.}}
(\byear{2013}).
\bhowpublished{Supplement to ``Limiting spectral distribution of sample
  autocovariance matrices.'' DOI:\doiurl{10.3150/13-BEJ520SUPP}}.
\bptok{imsref}%
\end{bmisc}
\endbibitem

%b8 ###
\bibitem{bosegangosen10}
\begin{barticle}[mr]
\bauthor{\bsnm{Bose},~\bfnm{Arup}\binits{A.}},
  \bauthor{\bsnm{Gangopadhyay},~\bfnm{Sreela}\binits{S.}} \AND
  \bauthor{\bsnm{Sen},~\bfnm{Arnab}\binits{A.}}
(\byear{2010}).
\btitle{Limiting spectral distribution of {$XX'$} matrices}.
\bjournal{Ann. Inst. Henri Poincar\'e Probab. Stat.}
\bvolume{46}
\bpages{677--707}.
\bid{doi={10.1214/09-AIHP329}, issn={0246-0203}, mr={2682263}}
\bptok{imsref}%
\end{barticle}
\endbibitem

%b9 ###
\bibitem{Bose08}
\begin{barticle}[mr]
\bauthor{\bsnm{Bose},~\bfnm{Arup}\binits{A.}} \AND
  \bauthor{\bsnm{Sen},~\bfnm{Arnab}\binits{A.}}
(\byear{2008}).
\btitle{Another look at the moment method for large dimensional random
  matrices}.
\bjournal{Electron. J. Probab.}
\bvolume{13}
\bpages{588--628}.
\bid{doi={10.1214/EJP.v13-501}, issn={1083-6489}, mr={2399292}}
\bptok{imsref}%
\end{barticle}
\endbibitem

%b10 ###
\bibitem{bottchersilberman}
\begin{bbook}[mr]
\bauthor{\bsnm{B{\"o}ttcher},~\bfnm{Albrecht}\binits{A.}} \AND
  \bauthor{\bsnm{Silbermann},~\bfnm{Bernd}\binits{B.}}
(\byear{1999}).
\btitle{Introduction to Large Truncated {T}oeplitz Matrices}.
\bseries{Universitext}.
\blocation{New York}: \bpublisher{Springer}.
\bid{doi={10.1007/978-1-4612-1426-7}, mr={1724795}}
\bptnote{check year}%
\bptok{imsref}%
\end{bbook}
\endbibitem

%b11 ###
\bibitem{bry}
\begin{barticle}[mr]
\bauthor{\bsnm{Bryc},~\bfnm{W{\l}odzimierz}\binits{W.}},
  \bauthor{\bsnm{Dembo},~\bfnm{Amir}\binits{A.}} \AND
  \bauthor{\bsnm{Jiang},~\bfnm{Tiefeng}\binits{T.}}
(\byear{2006}).
\btitle{Spectral measure of large random {H}ankel, {M}arkov and {T}oeplitz
  matrices}.
\bjournal{Ann. Probab.}
\bvolume{34}
\bpages{1--38}.
\bid{doi={10.1214/009117905000000495}, issn={0091-1798}, mr={2206341}}
\bptok{imsref}%
\end{barticle}
\endbibitem

%b12 ###
\bibitem{Dudley}
\begin{bbook}[mr]
\bauthor{\bsnm{Dudley},~\bfnm{R.~M.}\binits{R.M.}}
(\byear{2002}).
\btitle{Real Analysis and Probability}.
\bseries{Cambridge Studies in Advanced Mathematics}
\bvolume{74}.
\blocation{Cambridge}: \bpublisher{Cambridge Univ. Press}.
\bnote{Revised reprint of the 1989 original}.
\bid{doi={10.1017/CBO9780511755347}, mr={1932358}}
\bptok{imsref}%
\end{bbook}
\endbibitem

%b13 ###
\bibitem{GRS}
\begin{barticle}[mr]
\bauthor{\bsnm{Giraitis},~\bfnm{Liudas}\binits{L.}},
  \bauthor{\bsnm{Robinson},~\bfnm{Peter~M.}\binits{P.M.}} \AND
  \bauthor{\bsnm{Surgailis},~\bfnm{Donatas}\binits{D.}}
(\byear{2000}).
\btitle{A model for long memory conditional heteroscedasticity}.
\bjournal{Ann. Appl. Probab.}
\bvolume{10}
\bpages{1002--1024}.
\bid{doi={10.1214/aoap/1019487516}, issn={1050-5164}, mr={1789986}}
\bptok{imsref}%
\end{barticle}
\endbibitem

%b14 ###
\bibitem{hammil05}
\begin{barticle}[mr]
\bauthor{\bsnm{Hammond},~\bfnm{Christopher}\binits{C.}} \AND
  \bauthor{\bsnm{Miller},~\bfnm{Steven~J.}\binits{S.J.}}
(\byear{2005}).
\btitle{Distribution of eigenvalues for the ensemble of real symmetric
  {T}oeplitz matrices}.
\bjournal{J. Theoret. Probab.}
\bvolume{18}
\bpages{537--566}.
\bid{doi={10.1007/s10959-005-3518-5}, issn={0894-9840}, mr={2167641}}
\bptok{imsref}%
\end{barticle}
\endbibitem

%b15 ###
%  \bauthor{\bsnm{Pastur},~\bfnm{L.~A.}\binits{L.A.}}
%(\byear{1967}).

%b16 ###
\bibitem{mcmurraypolitis2010}
\begin{barticle}[mr]
\bauthor{\bsnm{McMurry},~\bfnm{Timothy~L.}\binits{T.L.}} \AND
  \bauthor{\bsnm{Politis},~\bfnm{Dimitris~N.}\binits{D.N.}}
(\byear{2010}).
\btitle{Banded and tapered estimates for autocovariance matrices and the linear
  process bootstrap}.
\bjournal{J. Time Series Anal.}
\bvolume{31}
\bpages{471--482}.
\bid{doi={10.1111/j.1467-9892.2010.00679.x}, issn={0143-9782}, mr={2732601}}
\bptok{imsref}%
\end{barticle}
\endbibitem

%b17 ###
\bibitem{SenA}
\begin{bmisc}[auto:STB|2013/06/05|13:45:01]
\bauthor{\bsnm{Sen},~\bfnm{Arnab}\binits{A.}}
(\byear{2006}).
\bhowpublished{Large dimensional random matrices. M. Stat. Project report, May
  2006. Indian Statistical Institute}.
\bptok{imsref}%
\end{bmisc}
\endbibitem

%b18 ###
\bibitem{SenS}
\begin{bmisc}[auto:STB|2013/06/05|13:45:01]
\bauthor{\bsnm{Sen},~\bfnm{Sanchayan}\binits{S.}}
(\byear{2010}).
\bhowpublished{Limiting spectral distribution of random matrices. M. Stat.
  Project report, July 2010. Indian Statistical Institute}.
\bptok{imsref}%
\end{bmisc}
\endbibitem

%b19 ###
\bibitem{WuPourahmadi}
\begin{barticle}[mr]
\bauthor{\bsnm{Wu},~\bfnm{Wei~Biao}\binits{W.B.}} \AND
  \bauthor{\bsnm{Pourahmadi},~\bfnm{Mohsen}\binits{M.}}
(\byear{2009}).
\btitle{Banding sample autocovariance matrices of stationary processes}.
\bjournal{Statist. Sinica}
\bvolume{19}
\bpages{1755--1768}.
\bid{issn={1017-0405}, mr={2589209}}
\bptok{imsref}%
\end{barticle}
\endbibitem

%b20 ###
\bibitem{xiaowu}
\begin{barticle}[mr]
\bauthor{\bsnm{Xiao},~\bfnm{Han}\binits{H.}} \AND
  \bauthor{\bsnm{Wu},~\bfnm{Wei~Biao}\binits{W.B.}}
(\byear{2012}).
\btitle{Covariance matrix estimation for stationary time series}.
\bjournal{Ann. Statist.}
\bvolume{40}
\bpages{466--493}.
\bid{doi={10.1214/11-AOS967}, issn={0090-5364}, mr={3014314}}
\bptnote{check year}%
\bptok{imsref}%
\end{barticle}
\endbibitem

%b21 ###
\bibitem{Yao2012}
\begin{barticle}[mr]
\bauthor{\bsnm{Yao},~\bfnm{Jianfeng}\binits{J.}}
(\byear{2012}).
\btitle{A note on a {M}ar\v cenko--{P}astur type theorem for time series}.
\bjournal{Statist. Probab. Lett.}
\bvolume{82}
\bpages{22--28}.
\bid{doi={10.1016/j.spl.2011.08.011}, issn={0167-7152}, mr={2863018}}
\bptok{imsref}%
\end{barticle}
\endbibitem

\end{thebibliography}
\end{document}